\documentclass{article}
\usepackage[top=1in, bottom=1in, left=1in, right=1in]{geometry}

\usepackage{amsfonts}
\usepackage{amsmath}
\usepackage{amssymb}
\usepackage{bm}
\usepackage{caption}
\usepackage{graphicx}
\usepackage{algorithm}
\usepackage{algpseudocode}
\usepackage{multirow}
\usepackage{bbm}
\usepackage{xcolor}
\usepackage{booktabs}
\usepackage{yfonts}
\usepackage{enumitem}
\usepackage{adjustbox}
\usepackage{mathtools}
\usepackage{float}
\usepackage{subfig}
\usepackage{verbatim}

\newcommand{\Be}{\bm{e}}

\newcommand{\Bp}{\bm{p}}

\newcommand{\R}{\ifmmode\mathbb{R}\else$\mathbb{R}$\fi}
\newcommand{\C}{\ifmmode\mathbb{C}\else$\mathbb{C}$\fi}
\newcommand{\N}{\ifmmode\mathbb{N}\else$\mathbb{N}$\fi}
\newcommand{\Q}{\ifmmode\mathbb{Q}\else$\mathbb{Q}$\fi}
\newcommand{\Z}{\ifmmode\mathbb{Z}\else$\mathbb{Z}$\fi}

\algnewcommand{\LeftComment}[1]{\Statex \(\triangleright\) #1}
\usepackage{hyperref}

\providecommand{\keywords}[1]
{
  \small	
  \textbf{\textit{Keywords---}} #1
}

\title{Solving High-Dimensional Partial Integral Differential Equations: The Finite Expression Method}

\author{
Gareth Hardwick{\thanks{Department of Mathematics, Purdue University}}, Senwei Liang{\thanks{Applied Mathematics and Computational Research Division, Lawrence Berkeley National Laboratory}} and Haizhao Yang{\thanks{Department of Mathematics and Department of Computer Science, University of Maryland College Park}}
}

\begin{document}
\maketitle
\noindent
\begin{abstract}
\noindent Partial integro-differential equations (PIDEs) have broad applications in the sciences, from electro-magnetism to options pricing.  In this paper, we introduce a new finite expression method (FEX) to solve PIDEs. This approach builds upon the original FEX and its inherent advantages with new advances: 1) A novel method of parameter grouping is proposed to reduce the number of coefficients in high-dimensional function approximation; 2) A Taylor series approximation method is implemented to significantly improve the computational efficiency and accuracy of the evaluation of the integral terms of PIDEs. The new FEX based method, denoted FEX-PG to indicate the addition of the parameter grouping (PG) step to the algorithm, provides both high accuracy and interpretable numerical solutions, with the outcome being an explicit equation that facilitates intuitive understanding of the underlying solution structures. These features are often absent in traditional methods, such as finite element methods (FEM) and finite difference methods, as well as in deep learning-based approaches. To benchmark our method against recent advances, we apply the new FEX-PG to solve benchmark PIDEs in the literature. In high-dimensional settings, FEX-PG exhibits strong and robust performance, achieving relative errors on the order of single precision machine epsilon, significantly outperforming existing approaches based on neural networks.
\end{abstract}


\section{Introduction}
Integro-differential equations, which involve both integrals and derivatives of a function have wide ranging applications from equations governing circuits \cite{lapa} to neuron behavior \cite{wilson}.  By extending these equations to functions of multiple variables, partial integro-differential functions (PIDEs) are introduced.  PIDEs have likewise proven critical across many applications in the sciences, from electro-magnetism \cite{kitora} to options pricing \cite{cruz}.  In lower dimensions, these PIDEs can be solved with mesh based approaches such as finite difference \cite{tang} and finite element methods.  However, high dimensional problems have proven much more difficult, as the mesh size (and general computational complexity) grows exponentially with the number of dimensions used.  Known as the ``curse of dimensionality", this phenomenon in its simplest form can be seen when estimating a general function to a given degree of accuracy - as the number of input variables increases, the number of samples required increases exponentially \cite{bellman}.  Lessening the ``curse" is particularly important in very high-dimensional problems, such as those encountered in finance, where each dimension might represent a different financial derivative in a portfolio \cite{doi:10.1073/pnas.1718942115}, and those in quantum chemistry, where the dimension is equal to the number of atoms in a system \cite{HAN2020109792}.  In such high-dimensional settings, traditional mesh-based methods become computationally infeasible, necessitating the exploration of alternative approaches. Deep learning has emerged as a powerful tool capable of solving problems in high-dimensional spaces \cite{Grohs_2023}, offering promising solutions where conventional techniques fall short.

Considerable progress has been made in the use of deep learning to solve differential equations~\cite{KHOO_LU_YING_2021, Li_2022,lu2021machinelearningellipticpdes,liang2022stiffness,liang2024solving}.  The use of mesh-free neural network (NN) based methods has proven very successful at solving PDEs of many kinds, in large part motivated by the super approximation power of NNs \cite{NEURIPS2020_2000f632, Shen_2021,Shen_2021_3layers,shen2022deepnetworkapproximationachieving}.  The Deep Ritz Method approximates the solution to a variational problem (such as partial differential equations, or PDEs) using an NN based on blocks of fully connected layers and residual connections to avoid the problem of vanishing gradients \cite{e2017deepritzmethoddeep, pmlr-v134-lu21a}.  Taking the Deep Ritz Method further, the Deep Nitsche Method changes the implementation of boundary conditions and demonstrated its ability to solve equations with 100 dimensions \cite{Ming_2021}.  The Deep Galerkin Method (DGM) swaps the linear combination of basis functions used by a traditional Galerkin Method with an NN used to approximate the solution.  By using random sampling points to train the NN, the need for a mesh is once again removed, allowing DGM to solve high dimensional PDEs \cite{Sirignano_2018}.  Another approach is physics informed informed machine learning \cite{Karniadakis2021}.  A physics informed neural network (PINN) implements the physical laws that govern the problem directly into the optimization loss function in the least square sense \cite{cuomo, raissi}.  Others have had success reformulating the PDEs as backward stochastic differential equations \cite{ E2017, doi:10.1073/pnas.1718942115}, using NNs to approximate the gradient of the solution.  While these approaches mitigate some of the problems arising from high dimensional spaces they introduce new challenges. The neural networks used by these approaches, especially when over-parameterized \cite{chizat}, can suffer from memory constraints and decreased accuracy at very high dimensions.  Additionally, the interpretability of solutions derived from neural networks remains limited, often presenting so-called ``black box" outcomes that obscure the underlying dynamics of the system.

To address the challenges associated with solving high-dimensional PDEs and PIDEs, the FEX method (Finite Expression) is introduced \cite{liang,song}.  The introduction of FEX is motivated by a few key characteristics: 1) FEX is very memory efficient; 2) FEX can achieve very high accuracy; 3) FEX gives a solution as a function written in standard notation with finitely many simple operators (i.e., a finite expression named in FEX), making the result interpretable to the user. A finite expression is a symbolic mathematical equation constructed using a finite number of operators, input variables, and constants.  The length of a finite expression is determined by the number of operators used, with a ``$k$-finite" expression having $k$ operators in total.  Previous research has demonstrated that these k-finite expressions are dense in various function spaces, ensuring that they can accurately represent solutions to complex problems \cite{liang,shen2022deepnetworkapproximationachieving}.  The operators used in these expressions, which include basic arithmetic operations and more complex functions (for example ``$+$", ``$-$", ``$\times$", ``$/$",``sin$(x)$" and ``$2^x$"), are categorized into binary and unary types. These operators are selected to generate an expression tree in computer algebra that formulates a mathematical expression.

The core of the FEX method involves selecting which of these operators to use, and where in the expression to place them in the expression tree.  In this way, the problem of identifying a finite expression to solve a PIDE is transformed into a combinatorial optimization (CO) problem \cite{liang} to select appropriate operators and constants.  Given this framework, the goal is to find the optimal sequence of operators that represents the solution to the equation.  Since the number of possible operator sequences is finite, this space can be searched systematically. This CO approach has already shown success in solving high-dimensional committor problems \cite{song}, general PDEs \cite{liang}, and in discovering physical laws from data \cite{jiang,song2024finiteexpressionmethodlearning}.

Given the demonstrated success of FEX, the goal to construct an FEX approach to solve PIDEs. The main challenge in developing the FEX lies scoring and optimizing the solutions.  To accomplish this, optimization must occur over two sets of variables: the set of operators chosen from when constructing the expression and the set of coefficients and constants multiplied and added to each term in the expression.  Thus, the FEX method operates in two stages: first, solving the combinatorial optimization problem to determine the correct sequence of operators, and second, learning the constants contained within the finite expression to achieve a high degree of accuracy in solving the PIDE.  This work seeks to explore the potential of FEX as a robust and interpretable tool for solving high-dimensional PIDEs, addressing the limitations of traditional methods and neural networks alike.

The adaptation of the FEX method to address high-dimensional PIDEs necessitated two significant advancements. The first is the development of an efficient technique for sampling the left- and right-hand sides of the PIDE, enabling the quantification of the accuracy of a proposed solution $u(t, x)$. Leveraging the structural properties of $u$ and the random variable of integration $z$, this approach utilizes a Taylor series expansion to approximate the integral term. The second advancement is the introduction of a novel method for dimensionality reduction in high-dimensional function approximations. This method substantially accelerates the learning process, both during the resolution of the combinatorial optimization problem to determine the correct operators and during the subsequent fine-tuning of the model parameters. The acceleration is achieved through the application of a modified linear layer, structured by a parameter hierarchy derived from a clustering algorithm.  Implementing these developments the new algorithm, FEX-PG is proposed, which enables the solution of the underlying equations with machine-epsilon accuracy, as demonstrated in this work.

\section{Preliminaries}
\subsection{Partial Integro-Differential Equation}
\label{sec:pide}

The primary aim of this paper is to find the solution $u(t, x)$ for $x \in \mathbb{R}^d$ and $0 < t < T$ to the following partial integro-differential equation~\cite{lu}: 
\begin{align}
\begin{cases} 
    \frac{\partial u}{\partial t} + b \cdot \nabla u + \frac{1}{2} \text{Tr}(\sigma \sigma^T H(u)) + \mathbf{A}u + f = 0, \\
    u(T, \cdot) = g(\cdot),
\end{cases}
\label{eqn:pide}
\end{align}
where $\nabla u$ and $H(u)$ represent the gradient and the Hessian matrix of $u(t, x)$ with respect to spatial variable $x \in \mathbb{R}^d$, respectively. Here $b = b(x) \in \mathbb{R}^d$ is a vector function, $\sigma = \sigma(x) \in \mathbb{R}^{d\times d}$ is a matrix function, and $f = f(t,x,u,\sigma^T \nabla u) \in \mathbb{R}$ is a given function. The operator $\mathbf{A}u$ is defined as:
\begin{align}
    \mathbf{A}u(t,x) = \int_{\mathbb{R}^d} \big(u(t,x+G(x,z)) - u(t,x) - G(x,z)\cdot \nabla u(t,x)\big)\nu(dz).
    \label{eqn:Au_integral}
\end{align}
Here $G = G(x,z) \in \mathbb{R}^d \times \mathbb{R}^{d} \rightarrow \mathbb{R}^d$, and $\nu$ is a L\'{e}vy measure associated with a Poisson random measure $N$.

For comparison, the NN method for solving Eqn.~\eqref{eqn:pide}, as implemented in \cite{lu}, is introduced, which will serve as a representative NN method for comparison in Section~\ref{sec:results}. In \cite{lu}, an NN is employed to approximate the solution  $u(x,t)$ with the trainable parameters optimized using a loss function derived from the forward-backward stochastic differential equations (FBSDE) system. 
 
 Specifically, the It\^{o}'s formula (as utilized in \cite{lu} and originally derived from \cite{Applebaum}) links the PIDE~\eqref{eqn:pide} to a corresponding set of stochastic differential equations. The following measures are introduced: 
\begin{align}
&N(t,S)(\omega) = \text{card}\{s\in [0,t):L_s(\omega)-L_{s-}(\omega) \in S\} \\
&\nu(S) = \mathbb{E}[N(1,S)(\omega)]\\
&\tilde{N}(t,S) = N(t,S)-t\nu(S)
\end{align}
$N(t,S)$ is the jump counting measure - it is exactly the number of jumps on the interval $(0,t]$ such that the jump size (given by the difference of the l\'{e}vy process $L_s$ and $L_{s-}$) is an element of the Borel set $S \in \mathbb{R}^d \setminus \{0\}$.  The jump measure $\nu$ is then used to create the compensated Poisson reandom measure $\tilde{N}$. Using these, It\^{o}'s formula  introduces the associated L\'{e}vy process:
 \begin{align}
&dX_t = b(X_t)dt +\sigma(X_t)dW_t + \int_{\mathbb{R}^d} G(X_t,z) \tilde{N}(dt,dz),\\
&Y_t = u(t, X_t),\\
&Z_t = \nabla u(t, X_t),\\
&U_t = \int_{\mathbb{R}^d}(u(t, X_t + G(X_t, z)) - u(t, X_t)))\nu(dz).
 \end{align}
   Hence, FBSDE is given by:
\begin{align}
&dX_t = b(X_t)dt +\sigma(X_t)dW_t + \int_{\mathbb{R}^d} G(X_t,z) \tilde{N}(dt,dz),\\
\begin{split}
&dY_t = -f(t, X_t, Y_t, \sigma(X_t)^T Z_t)dt + (\sigma(X_t)^T Z_t)^TdW_t \\
&+ \int_{\mathbb{R}^d}[u(t, X_t + G(X_t, z)) - u(t, X_t)]\tilde{N}(dt, dz).
\label{eqn:dy}
\end{split}
\end{align}
The fundamental connection between this system and the PIDE~\eqref{eqn:pide} is that if $u(t,x)$ satisfies the PIDE, then the processes $(X_t, Y_t, Z_t, U_t)$ must satisfy the FBSDE. With an NN function $\mathrm{NN}(x,t)$ approximating $u(x,t)$, trajectories of the FBSDE can be simulated with the approximate solution $\mathrm{NN}(x,t) \approx u(x,t)$. The difference between $Y_{t_{n}} \approx u(X_{t_{n}}, t_{n})$ and $Y_{t_{n+1}} \approx u(X_{t_{n+1}}, t_{n+1})$ should align with Eqn.~\eqref{eqn:dy}, which establishes a temporal difference loss function to optimize the NN solution in~\cite{lu}.

\subsection{Finite Expression Method}
\label{sec:fex}
\begin{figure}[!ht]
\includegraphics[scale = .6]{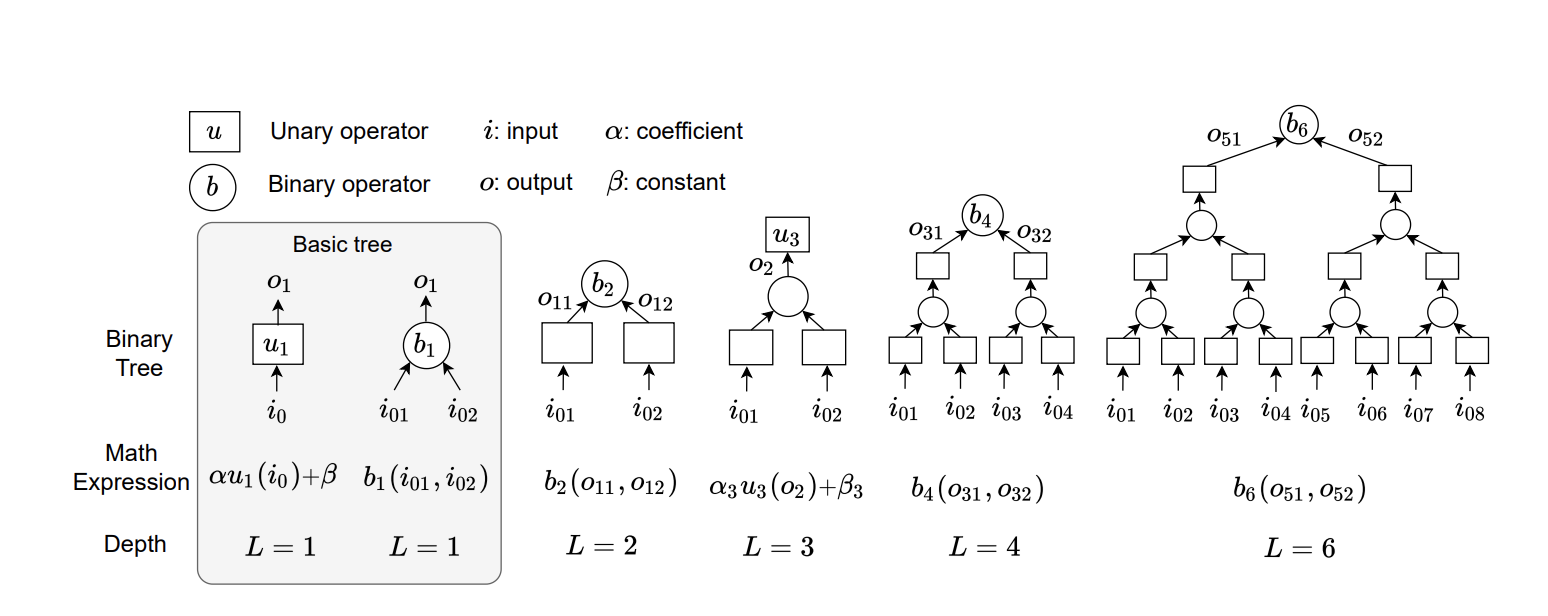}
\caption{ Computation structure using binary trees.  Each node is either a binary or unary operator.  Beginning with depth-1 trees, mathematical expressions can be built by performing computation recursively. Each tree node is either a binary operator or a unary operator that takes value from the corresponding binary or unary set. The binary set can be $\mathbb{B}:=\{+,-,\times,\div,\cdots\}$. The unary set can be $\mathbb{U}:=\{\sin,\exp, \log, \text{Id}, (\cdot)^2, \int\cdot\text{d} x_i, \frac{\partial\cdot}{\partial x_i}, \cdots\}$, which contains elementary functions (e.g., polynomial and trigonometric function), antiderivative and differentiation operators. Here ``Id'' denotes the identity map. Notice that if an integration or a derivative is used in the expression, the operator can be applied numerically. Reproduced from \cite{liang}, with permission.}
\label{fig:tree}
\end{figure}

FEX seeks a solution to a PDE in the function space of mathematical expressions composed of a finite number of operators. In FEX implementation, a finite expression is represented as a binary tree $\mathcal{T}$, as shown in Figure~\ref{fig:tree}. Each node in the tree is assigned a value from a set of operators, forming an operator sequence $\Be$. Each operator is associated with trainable scaling and bias parameters, denoted by $\bm{\theta}$. Thus, a finite expression can be represented as $u(\bm{x}; \mathcal{T}, \mathbf{e}, \bm{\theta})$. The goal is to identify the mathematical expression by minimizing the functional $\mathcal{L}$ related to a PDE, where the minimizer of $\mathcal{L}$ corresponds to the solution of the PDE. Specifically, the resulting combinatorial optimization problem is:
\begin{align}
\min \{\mathcal{L}(u(\cdot; \mathcal{T}, \Be, \bm{\theta}))|\Be, \bm{\theta}\}.
\label{eqn:obj}
\end{align}

\begin{figure}[!ht]
\begin{center}
\includegraphics[scale = .55]{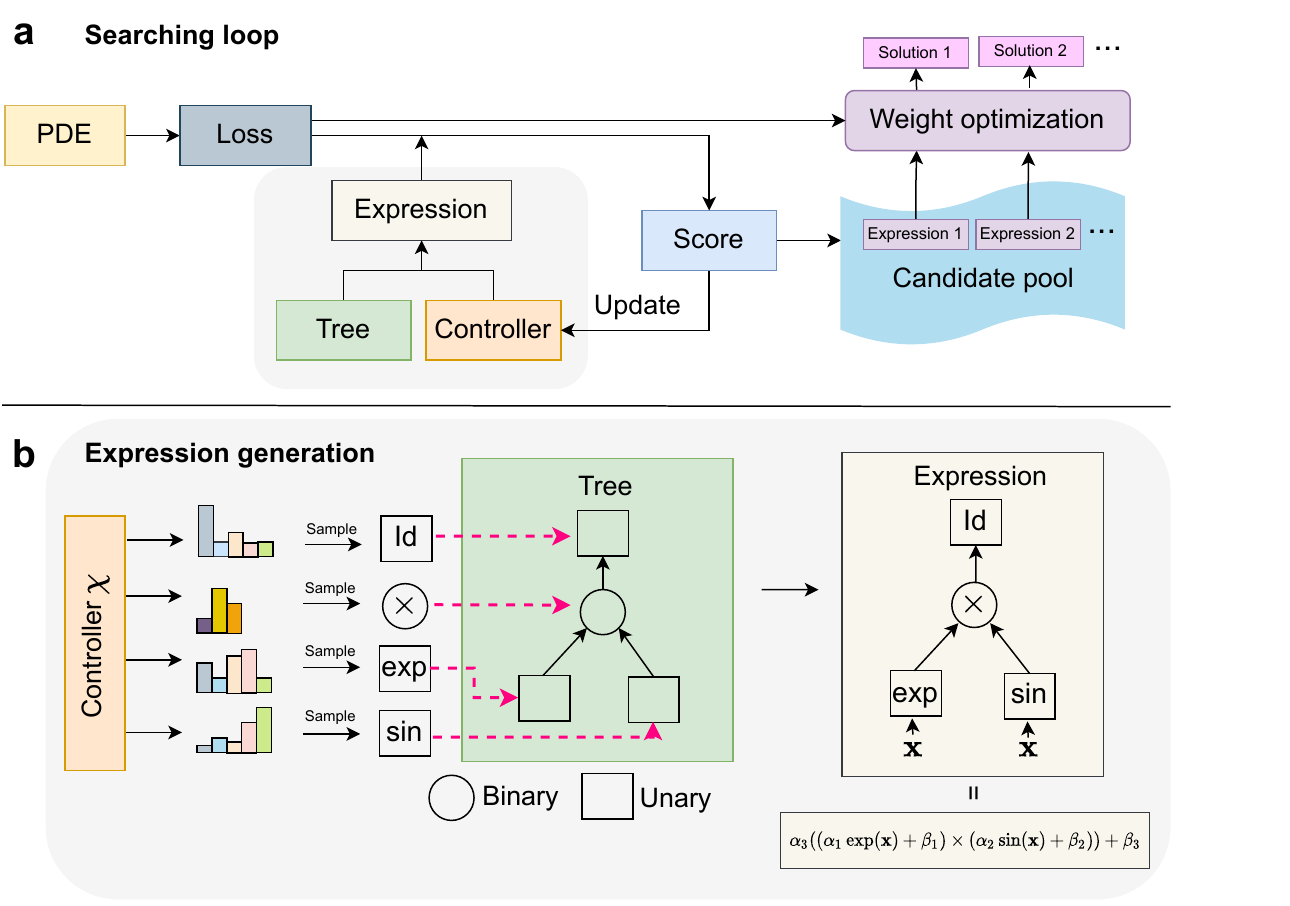}
\caption{A flowchart outlining the FEX algorithm: (a) The search loop consists of four key components: score computation, operator sequence generation, controller updates, and candidate optimization. (b) A schematic illustration of the controller used to generate an operator sequence for a mathematical expression. Reproduced from \cite{liang}, with permission.}
\label{fig:FEXdiagram}
\end{center}
\end{figure}

In FEX, to address this CO, a search loop (see Fig.~\ref{fig:FEXdiagram}a) based on reinforcement learning is employed to identify effective operators $\Be$ that can potentially recover the true solution when selected in the expression. In FEX, the search loop consists of four main components:

\begin{enumerate}
    \item \textbf{Score computation (i.e., rewards in RL)}: A mixed-order optimization algorithm is introduced to efficiently evaluate the score of the operator sequence $\Be$, helping to reveal the true structure. A higher score indicates a greater likelihood of identifying the true solution.

The score of $\Be$, denoted as $S(\Be)$, is defined on the interval $[0,1]$ by:
\begin{align}
S(\Be) := \big(1+L(\Be)\big)^{-1},
\label{eqn:orgscore}
\end{align}
where $L(\Be) := \min \{\mathcal{L}(u(\cdot; \mathcal{T}, \Be, \bm{\theta}))|\bm{\theta}\}$. As $L(\Be)$ approaches 0, the expression represented by $\Be$ comes closer to the true solution, causing the score $S(\Be)$ to approach 1. Conversely, as $L(\Be)$ increases, $S(\Be)$ approaches 0.
Finding the global minimizer of $\mathcal{L}(u(\cdot; \mathcal{T}, \Be, \bm{\theta}))$ with respect to $\bm{\theta}$ is computationally expensive and challenging. To expedite the evaluation of $S(\Be)$, rather than conducting an exhaustive search for a global minimizer, FEX employs a combination of first-order and second-order optimization algorithms. The optimization process consists of two stages. First, a first-order algorithm is employed for $T_1$ iterations to obtain a well-informed initial estimate. This is followed by a second-order algorithm (such as BFGS~\cite{fletcher2013practical}) for an additional $T_2$ iterations to further refine the solution.
Let $\bm{\theta}_0^{\Be}$ denote the initial parameter set, and $\bm{\theta}_{T_1+T_2}^{\Be}$ represent the parameter set after completing $T_1+T_2$ iterations of this two-stage optimization process. The resulting $\bm{\theta}_{T_1+T_2}^{\Be}$ serves as an approximation of $\arg \min_{\bm{\theta}} \mathcal{L}(u(\cdot; \mathcal{T}, \Be, \bm{\theta}))$.
Then, $S(\Be)$ is estimated by:
\begin{align}
S(\Be) \approx \big(1+\mathcal{L} (u(\cdot; \mathcal{T}, \Be, \bm{\theta}_{T_1+T_2}^{\Be}))\big)^{-1}.
\label{eqn:score}
\end{align}

 \item \textbf{Operator sequence generation (i.e., taking actions in RL)}:  The controller is to generate high-scoring operator sequences during the search process (see Fig.~\ref{fig:FEXdiagram}b). We denote the controller as $\bm{\chi}_\Phi$, where $\Phi$ represents its model parameters. Throughout the search, $\Phi$ is updated to increase the likelihood of producing favorable operator sequences. The process of sampling an operator sequence $\Be$ from the controller $\bm{\chi}_\Phi$ is denoted as $\Be\sim\bm{\chi}_\Phi$.
Considering the tree node values of $\mathcal{T}$ as random variables, the controller $\bm{\chi}_\Phi$ outputs a series of probability mass functions $\Bp_\Phi^1, \Bp_\Phi^2, \cdots, \Bp_\Phi^s$ to characterize their distributions, where $s$ represents the total number of nodes. Each tree node value $e_j$ is sampled from its corresponding $\Bp_\Phi^j$ to generate an operator. The resulting operator sequence $\Be$ is then defined as $(e_1, e_2, \cdots, e_s)$.
To enhance the exploration of potentially high-scoring sequences, an $\epsilon$-greedy strategy is used. With a probability of $\epsilon < 1$, $e_i$ is sampled from a uniform distribution over the operator set. Conversely, with a probability of $1-\epsilon$, $e_i$ is sampled from $\Bp_\Phi^i$. A higher value of $\epsilon$ increases the likelihood of exploring new sequences.

 \item \textbf{Controller update (i.e., policy optimization in RL)}: The controller is updated to increase the probability of generating better operator sequences based on the feedback from their scores. While various methods (e.g., heuristic algorithms) can be used to model the controller, the policy gradient method from RL is employed for optimization.

FEX adopts the objective function proposed by \cite{petersen2021deep} to update the controller. This function is 
\begin{align}
\mathcal{J}(\Phi)=\mathbb{E}_{\Be \sim \bm{\chi}_\Phi} \{S(\Be)|S(\Be)\geq S_{\nu, \Phi}\},
\label{eqn:expectriskseeking}
\end{align}
where $S_{\nu, \Phi}$ denotes the $(1-\nu)\times 100\%$-quantile of the score distribution generated by $\bm{\chi}_\Phi$ within a given batch. This objective function focuses on optimizing the upper tail of the score distribution, thereby increasing the likelihood of discovering high-performing operator sequences.

The controller parameter $\Phi$ is updated via the gradient ascent with a learning rate $\eta$, i.e., 
\begin{align}
\Phi \leftarrow \Phi+\eta \nabla_\Phi\mathcal{J}(\Phi). 
\label{eqn:gradientascent}
\end{align}

 \item \textbf{Candidate optimization (i.e., a policy deployment)}: During the search, a candidate pool is maintained to store high-scoring operator sequences. After the search, the parameters $\bm{\theta}$ of the high-scoring operator sequence $\Be$ are optimized to approximate the PDE solution.

The score of an operator sequence $\Be$ is determined by optimizing a nonconvex function, starting from a random initial point and using a limited number of update iterations. This approach has limitations: the optimization process may become trapped in suboptimal local minima, and consequently, the score may not always accurately reflect how well $\Be$ captures the structure of the true solution. In fact, the operator sequence that closely approximates or exactly matches the true solution may not necessarily achieve the highest score.
To mitigate the risk of overlooking promising operator sequences, we maintain a candidate pool $\mathbb{P}$ with a fixed capacity $K$. This pool is designed to store multiple high-scoring sequences of $\Be$.

After the search loop concludes, we perform an additional optimization step for each $\Be$ in $\mathbb{P}$. Specifically, we optimize the objective function $\mathcal{L}(u(\cdot; \mathcal{T}, \Be, \bm{\theta}))$ with respect to $\bm{\theta}$ using a first-order algorithm. This optimization runs for $T_3$ iterations with a small learning rate.

\end{enumerate}


\section{Proposed FEX-PG Method for PIDEs}
\label{sec:alg}
This section introduces the proposed functional $\mathcal{L}$~\eqref{eqn:obj} used in FEX-PG to identify the solution to the PIDE~\eqref{eqn:pide}. Additionally, a method for efficient evaluation of the integral terms is presented. Finally, to improve efficiency in FEX, a method for grouping certain trainable parameters is introduced to simplify the problem. The proposed FEX-PG algorithm is summarized in Alg. \ref{alg1}.

\subsection{Functional for the PIDE Solution}
To apply the FEX to solve PIDEs, we propose a functional used to evaluate a candidate function. This functional consists of a least squares loss which combines both the equation loss and boundary loss.

Denote Eqn.~\eqref{eqn:pide} as $\mathcal{D}(u)=0$. The equation loss is defined as $\|\mathcal{D}(u)\|_{L_2([0,T] \times \Omega)}$. To enforce the boundary condition, the boundary loss is defined as $\|u(T, \cdot) - g(\cdot)\|_{L_2(\Omega)}$. Thus, the function $\mathcal{L}$ is given by the following:
\begin{align*}
    \mathcal{L}(u)= \|\mathcal{D}(u)\|^2_{L_2([0,T]\times \Omega)} + \|u(T, \cdot)-g(\cdot)\|^2_{L_2(\Omega)}.
\end{align*}
Such a loss can be approximated using $N$ random points $(t_i, x_i)$ within the domain, where $t_i \in [0,T]$ and $x_i \in \Omega$, and $M$ points $(T, x_j)$ on the boundary with $x_j \in \Omega$. Hence the functional is approximated by:
\begin{align}
    \mathcal{L}(u) \approx \frac{1}{N}\sum_{i=1}^N|\mathcal{D}(\tilde{u}(t_i,x_i))|^2+ \frac{1}{M}\sum_{j=1}^M | \tilde{u}(T, x_j) - g(x_j)|^2.
    \label{eqn:leastsquare}
\end{align}
Once $\mathcal{L}$ is defined, it can be employed with FEX, as introduced in Section~\ref{sec:fex}, to find the solution. 



\subsection{Evaluation of the Integral Term}
\label{sec:evalint}
To compute the integral term in Eqn.~\eqref{eqn:Au_integral} within the proposed least squares loss~\eqref{eqn:leastsquare} for each sample point $(t_i, x_i)$, efficient evaluation of the integral is crucial. This section will therefore present an estimation for simplifying the calculation of the integral term.

Starting with some manipulation of the integral:
\begin{align}
    \mathbf{A}u(t,x)&=\int_{\mathrm{R}^d} (u(t, x+ G(x,z)) - u(t,x) - G(x,z) \cdot \nabla  {u}(t,x))\nu(dz) \\
    & = \int_{\mathrm{R}^d}  {u}(t, x+G(x,z)) -  {u}(t,x) \nu(dz) - \int_{\mathrm{R}^d} G(x,z) \cdot \nabla  {u}(t,x) \nu(dz) \\
    & = \lambda \int_{\mathrm{R}^d} ( {u}(t, x+G(x,z)) -  {u}(t,x))\phi(z)dz - \lambda \nabla  {u}(t,x) \cdot \int_{\mathrm{R}^d} G(x,z) \phi(z)dz \\
    & = \lambda \Big(\mathbb{E}[ {u}(t,x+G(x,z))] -  {u}(t,x) - \mathbb{E}[G(x,z)] \cdot \nabla  {u}(t,x) \Big),
    \label{eqn:integral}
\end{align}
where the expectation is taken with respect to the random variable $z$. 


In \eqref{eqn:integral}, evaluating $\mathbb{E}[{u}(t, x + G(x,z))]$ is particularly challenging, especially in high-dimensional cases. In lower dimensions, this can be handled using standard quadrature methods, such as the trapezoidal rule. However, as the dimensionality increases, these methods become impractical. Even a mesh free Monte Carlo evaluation is likewise prohibitively expensive in high dimensions as the number of required sampling points still increases exponentially with dimension. To address this, the expected value of the Taylor series expansion of the candidate function, centered at $(t, x+\mathbb{E}[G(x,z)])$ is instead computed \cite{OLSON1991309}. In the case of Eqn.~\eqref{eqn:3.9}, when $G(x,z) = z$, and all dimensions of $z\in \mathrm{R}^d$ have identical mean and variance $(\mu, \sigma^2)$, based on Taylor expansion, 
\begin{align}
    {u}(t, x + z) = {u}(t,x+\mu) + \nabla {u}(t,x+\mu)^T (z-\mu)
    + \frac{1}{2} (z-\mu)^T \textbf{H}(u)(z-\mu) + \cdots. 
\end{align}
The expectation with regards to the space variable $z$ is taken, resulting in
\begin{align}
\mathbb{E}[{u}(t, x + z)] \approx {u}(t,x+\mu) + \frac{1}{2}\sigma^2 \sum_{i=1}^d \textbf{H}(u)_{ii}(t,x+\mu).
\label{eqn:intgralest}
\end{align}
Note that the term involving the gradient disappears when we take the expectation as $\mathbb{E}[z - \mu] = \mathbb{E}[z - \mathbb{E}[z]] = 0$ by linearity of the expectation. With this estimation, evaluating $\mathbb{E}[{u}(t, x + G(x,z))]$ is much more efficient. A bound on the error from this approximation has been established for specific functions $u$ \cite{sym12081379, vershynin2011introductionnonasymptoticanalysisrandom}, but further work is required to establish a bound given a general candidate ${u}$.  However, numerical results indicate that the error remains within an acceptable range. While it is possible to include even higher-order terms, no improvement in performance was found, and the computational time increased significantly when computing $4^{th}$-order derivatives.  Note that as $\sigma$ grows, so do the higher order terms that are not included in this truncated form.  Specifically, while the odd central moments of the normal distribution are all zero, the $2n^{th}$ terms grow with $\mathcal{O}(\sigma^{2n})$.  Clearly, if $\sigma \geq 1$, this may prove to be a large source of inaccuracy if the higher order derivatives of the candidate function cannot control the these terms.  However, as seen in Section~\ref{sec:highdpide}, numerical stability of the proposed FEX-PG (Alg. \ref{alg1}), up to $\sigma = 1$ is observed.

\subsection{Parameter Grouping}
\label{sec:paramgrouping}
To enhance the efficiency of the search loop in FEX, a method is proposed for grouping certain scaling parameters in the leaf nodes into clusters, with parameters within each group being shared.

Each of the leaves of the binary tree structure takes as input $(t, x)$, where $x:=[x_1, \cdots, x_d] \in \mathbb{R}^d$, and applies a unary function element-wise to each input along with some learnable parameters.  Letting the unary function of the $i^{th}$ leaf be $\phi_i$, and its weights be $\alpha_i:=[\alpha_i^0, \cdots, \alpha_i^d]\in \mathbb{R}^{d+1}$ and $\beta_i\in \mathbb{R}$, the output of the leaf can be written as 
\begin{align}
   \alpha_i^0 \phi_i(t) + \alpha_i^1 \phi_i(x_0) + \alpha_i^2 \phi_i(x_1) + \cdots + \alpha_i^{d} \phi_i(x_d).
\end{align}

In practice, some of the $\alpha$ values may be similar or potentially identical. It is possible to leverage this to minimize the number of unique coefficients, which simplifies and accelerates the learning process. To achieve this, the weight vector $\alpha_i$ is input to an unsupervised clustering algorithm. The output is a vector indicating the cluster assignment for each element of $\alpha_i$.
\begin{figure}[!ht]
\begin{center}
\includegraphics[scale=.2]{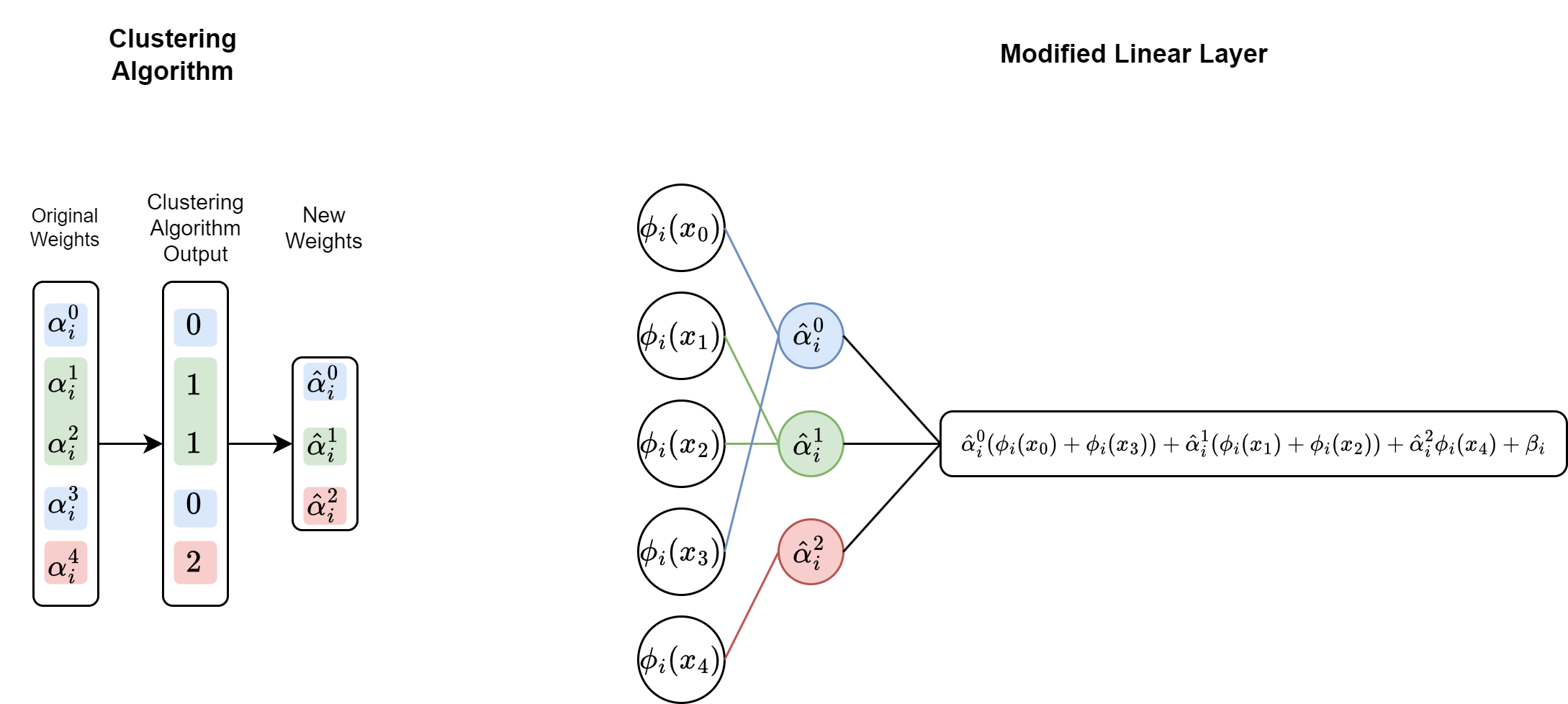}
\caption{Example of parameter grouping. In this example, the clustering algorithm groups $\alpha_i^0$ with $\alpha_i^3$ and $\alpha_i^1$ with $\alpha_i^2$, while $\alpha_i^4$ remains unclustered (in this case, three coefficients are needed)}
\label{fig:group}
\end{center}
\end{figure}
Figure~\ref{fig:group} illustrates the application of clustering to the weights of the $i^{th}$ leaf in the tree. The clustering algorithm identifies which input dimensions should share the same coefficient, and the total number of unique coefficients can be determined by adding one to the largest value in the output vector. A new linear layer is constructed that follows the structure provided by the clustering algorithm's output (as shown in Figure~\ref{fig:group} (Right)). After this clustering is performed, re-learning the weights and the bias term $\beta_i$ from scratch ensures that no biases is carried over from the previous formulation of the solution. In practice, SciPy's hierarchical clustering algorithm \cite{müllner2011modernhierarchicalagglomerativeclustering, 2020SciPy-NMeth} is used, with a threshold parameter set to $\eta:=\frac{1}{d^2}$.  The threshold parameter determines the maximum 
{ cophenetic distance
}
between elements in a cluster, in practice this hyperparameter is problem dependent. 
\begin{figure}[H]
\begin{center}
\includegraphics[scale = .48]{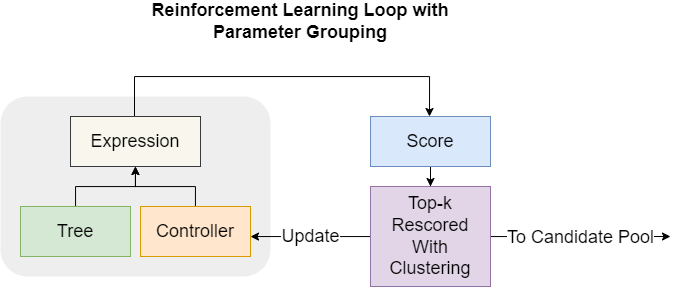}
\caption{Searching loop with parameter grouping. }
\label{fig:searchinggroup}
\end{center}
\end{figure}
By grouping the parameters, structure is introduced (through parentheses that distribute a single weight across multiple input dimensions, as in Figure~\ref{fig:group}) and this reduces the dimensionality of the optimization problem during fine-tuning. With fewer parameters, the process of fitting constants and coefficients becomes significantly faster. To fully capitalize on this advantage, this parameter grouping is incorporated into the search loop itself. Figure \ref{fig:searchinggroup} illustrates how this grouping is applied within the reinforcement learning loop for searching. The application within FEX-PG is straightforward: during the loop, the expressions are first scored as usual. Then, the top-scoring expression (denoted as $\Be^{\ast}$) is taken, and the restructuring process is applied. Following that, $T_3$ iterations of Adam are performed to fine-tune the restructured tree's weights.  In practice, $T_3 = 100$ iterations were found to be more than sufficient to compute a new, more accurate score. The original loss for $\Be^{\ast}$ is then replaced with the newly optimized loss from the tuned tree. Importantly, the output of the clustering algorithm is saved
 to the candidate pool. This ensures that the computationally expensive restructuring process does not need to be repeated during fine-tuning.


 The advantage of this implementation is that the sequence $\Be^{\ast}$ is scored much more accurately, which significantly improves the quality of the update to the controller $\chi$, allowing the searching process to become much more efficient. In the worst-case scenario, where the restructuring process fails and the new loss is worse than the previous one, the updated loss is discarded and algorithm proceeds as normal, ensuring no negative impact on the overall process.
\begin{algorithm}
\begin{algorithmic}[1]
\caption{Fixed Tree FEX-PG for PIDEs}
\LeftComment{Input: PIDE; A tree $\mathcal{T}$; Searching loop iteration $T$; Coarse-tune iteration $T_1$ with Adam; Coarse-tune iteration $T_2$ with BFGS; Medium-tune iteration $T_3$ with Adam; Fine-tune iteration $T_4$ with Adam; Pool size $K$; Batch size $N$; Clustering threshold $\eta$.}
\LeftComment{Output: The solution $u(\textbf{x}; \mathcal{T}, \hat{\Be}, \hat{\theta})$}
\State Initialize the agent $\chi$ for the tree $\mathcal{T}$
\State $\mathbb{P} \leftarrow \{\}$
\For{$\hbox to 1em{\thinspace\hrulefill\thinspace}$ from 1 to $T$}
    \State Sample $N$ sequences $\{\Be^{(1)}, \Be^{(2)},...,\Be^{(N)}\}$ from $\chi$
    \State Losses $\leftarrow [\text{ }]$
    \For{n from 1 to $N$}
        \State Minimize $\mathcal{L}(u(\cdot; \mathcal{T}, \Be^{(n)}, \theta^{(n)}))$ with respect to $\theta^{(n)}$ by coarse-tune with $T_1 + T_2$ iterations
        \State After $T_1 + T_2$ iterations, Losses.append($\mathcal{L}(u(\cdot; \mathcal{T}, \Be^{(n)}, \theta_{T_1 + T_2}^{(n)}))$)
    \EndFor
    \State Denote $ \Tilde{n} := \arg\min (\text{Losses})$
    \State Apply operator sequence $\Be^{(\Tilde{n})}$ to tree $\mathcal{T}$, denoted as $\mathcal{T}_{e^{(\Tilde{n})}}$
    \For{leaf in $\mathcal{T}_{e^{(\Tilde{n})}}$} \Comment{Parameter Grouping}
        \State  Apply hierarchical clustering algorithm with threshold parameter $\eta$
        \State Replace the linear layer of each leaf with the modified linear layer (as depicted in Figure~\ref{fig:group})
    \EndFor
    \For{$\hbox to 1em{\thinspace\hrulefill\thinspace}$  from 1 to $T_3$} \Comment{Learning weights for new modified linear layers}
        \State Calculate $ \mathcal{L}(u(\cdot;\mathcal{T}_{e^{(\Tilde{n})}}, e^{(\Tilde{n})}, \theta^{(\tilde{n})}))$ using $\mathcal{T}_{e^{(\Tilde{n})}}$ and update $\theta$ with Adam
        \If{$\hbox to 1em{\thinspace\hrulefill\thinspace} = T_3$ and Losses[$\tilde{n}$] $<$ $\mathcal{L}(u(\cdot;\mathcal{T}_{e^{(\Tilde{n})}}, e^{(\Tilde{n})}, \theta_{T_3}^{(\tilde{n})}))$}
            \State $\text{Losses}[\Tilde{n}] \leftarrow  \mathcal{L}(u(\cdot;\mathcal{T}_{e^{(\Tilde{n})}}, e^{(\Tilde{n})}, \theta_{T_3}^{(\tilde{n})}))$
            \EndIf
    \EndFor
    \State Calculate rewards using Losses[:] and update $\chi$
    \For{n from 1 to $N$}        
        \If{Losses$[n] <$ any in $\mathbb{P}$}
            \State $\mathbb{P}$.append($\Be^{(n)}$)
            \State $\mathbb{P}$ pops $\Be$ with the smallest reward when overloading
            \EndIf
    \EndFor
\EndFor
\For{$\Be$ in $\mathbb{P}$} \Comment{Candidate optimization}
    \For{$\hbox to 1em{\thinspace\hrulefill\thinspace}$  from 1 to $T_4$}
        \State Minimize $\mathcal{L}(u(\cdot; \mathcal{T}, \Be, \theta))$ with respect to $\theta$ using Adam
        \If{all(previous 5 values of $\mathcal{L}$) $< 1e-14$}
            \State Break
        \EndIf
    \EndFor
\EndFor
\State \textbf{Return} the expression with the smallest fine-tune error
\label{alg1}
\end{algorithmic}
\end{algorithm}
\section{Numerical Results}
\label{sec:results}
 This section presents numerical verification of the effectiveness of the proposed methods for solving the PIDE~\eqref{eqn:pide}. First, two 1-dimensional cases are tested in Section~\ref{sec:1d}, followed by two high-dimensional cases (up to 100 dimensions) in Section~\ref{sec:highdpide}. Since traditional methods are prone to the curse of dimensionality, the proposed FEX-PG method is compared to with the TD-based neural network (referred as TD-NN) method \cite{lu} as introduced in Section~\ref{sec:pide}. The test cases used in ref. \cite{lu} are adopted for a true apples-to-apples comparison. In all of these examples, the same set of binary operators: ``$+$", ``$-$" and ``$\times$", and unary operators: ``$0$", ``$1$", ``$x$", ``$x^2$'', ``$x^3$", ``$x^4$", ``$e^x$", ``$\sin (x)$" and ``$\cos (x)$" are used.  Along with these operators, a depth-2 tree structure is used, which can be seen on the far right of Figure~\ref{fig:tree}. 

\subsection{Low Dimensional PIDE}
\label{sec:1d}
We begin the discussion by referencing two one-dimensional PIDEs from \cite{lu}. Note that in this case, there is no need to apply parameter grouping or the integral estimation discussed in Section~\ref{sec:alg}. This aims to demonstrate that the vanilla FEX method is capable of solving PIDEs for one dimensional case. {
For the integral term, rather than applying the integral estimation from Section~\ref{sec:evalint}, a trapezoidal rule is used to calculate the integral.  No stopping criteria other than a maximum number of iterations was used for these examples as they serve to demonstrate vanilla FEX's potential at solving PIDEs.
}

\textbf{Example 1. } The first PIDE is given by
\begin{align} \begin{cases} 
    \frac{\partial u}{\partial t} + \int_{\mathrm{R}} \big(u(t, xe^z) - u(t,x) - x(e^z-1) \frac{\partial u}{\partial x}\big)\nu(dz) = 0,\\
    u(T, x) = x.\\
    \end{cases}
    \label{eqn:1d_1}
\end{align}
Here, $\nu(dz) = \lambda\phi(z)dz$, with $\phi(z) = \frac{1}{\sqrt{2 \pi}\sigma}e^{-\frac{1}{2}\left(\frac{z-\mu}{\sigma}\right)^2}$ and $x \in \mathbb{R}$. The same parameters as in \cite{lu} are used. The intensity of the Poisson process is set to $\lambda = 0.3$, and the jumps follow a normal distribution with $\mu = 0.4$ and $\sigma = 0.25$.  Both $x \text{ and } t \in [0,1]$, with $T = 1$. 
{
In practice, it was found that using 50 evenly spaced points in $[0,1]$ for numerical integration (via the trapezoid rule) was sufficient.
}
The true solution of Eqn. \eqref{eqn:1d_1} is $u(t,x) = x$.
\begin{figure}[H]
    \centering
    \subfloat[\centering Loss During Fine-tuning]{{\includegraphics[width=7cm]{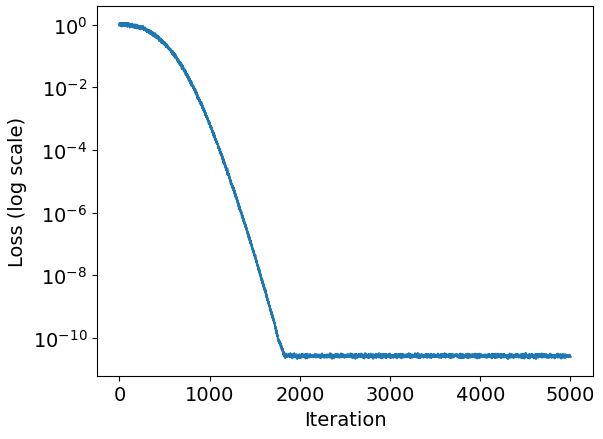} }}%
    \qquad
    \subfloat[\centering Relative Error During Fine-Tuning]{{\includegraphics[width=7.1cm]{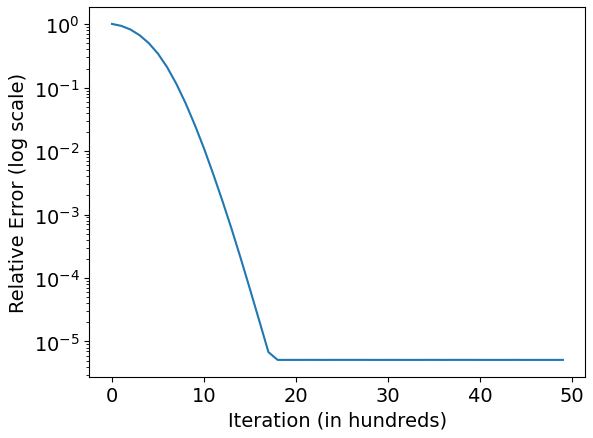}}}%
    \caption{Optimization profile for PIDE~\eqref{eqn:1d_1}. \textbf{(a):} Training loss of the candidate solution during fine-tuning. \textbf{(b):} Relative error of the candidate solution. }%
    \label{fig:3.2loss}%
\end{figure}


After just 50 iterations of the search loop, promising candidate functions are identified. Figure \ref{fig:3.2loss}(a) shows the loss of the candidate function as it is fine-tuned over iterations, while Figure \ref{fig:3.2loss}(b) displays the corresponding relative error. It is evident that after only 2000 iterations of fine-tuning, the relative error drops to the order of $10^{-6}$. The final solution obtained is $\tilde{u}(t,x) = 0.9999997x + 0.0000837$, which is very close to the true solution $u(t,x) = x$.

\textbf{Example 2. } Next, consider another example from \cite{lu}. The PIDE is given by: 
\begin{align} \begin{cases} 
    \frac{\partial u}{\partial t} + \epsilon x  \frac{\partial u}{\partial x} + \frac{1}{2}\theta^2 \frac{\partial^2 u}{\partial x^2} + \int_{\mathrm{R}} \big(u(t, xe^z) - u(t,x) - x(e^z-1) \frac{\partial u}{\partial x}\big)\nu(dz) = \epsilon x,\\
    u(T, x) = x.
    \label{eqn:3.5}
    \end{cases}
\end{align}
Here, $\nu(dz) = \lambda\phi(z)dz$, with $\phi(z) = \frac{1}{\sqrt{2 \pi}\sigma}e^{-\frac{1}{2}\left(\frac{z-\mu}{\sigma}\right)^2}$ and $z \in \mathbb{R}$. Again, the same parameters and domain as in \cite{lu} are used: $\lambda = 0.3$, $\mu = 0.4$, $\sigma = 0.25$, $\epsilon = 0.25$ and $\theta = 0$ with  $x, t \in [0,1]$, and $T=1$.  
{
Once again, 50 points evenly spaced in $[0,1]$ were sufficient for numerical integration.  
}
The true solution is $u(t,x) = x$.

\begin{figure}[H]
    \centering
    \subfloat[\centering Loss During Fine-tuning]{{\includegraphics[width=7cm]{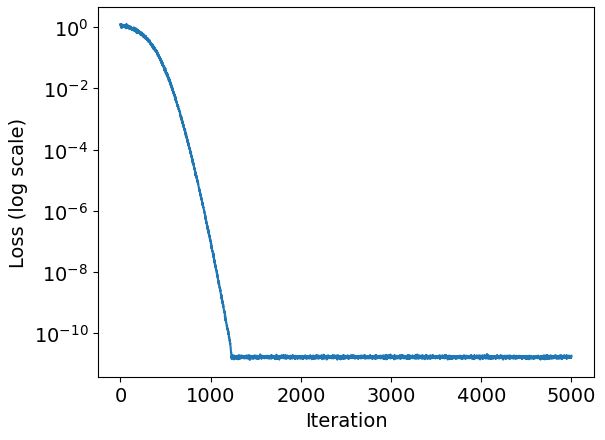} }}%
    \qquad
    \subfloat[\centering Relative Error During Fine-tuning]{{\includegraphics[width=7cm]{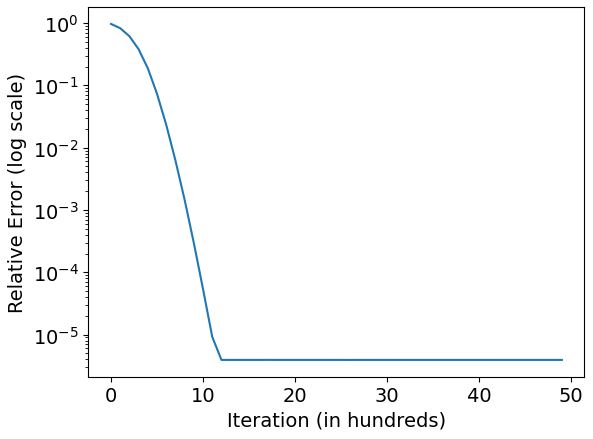} }}%
    \caption{Optimization profile for PIDE~\eqref{eqn:3.5}. \textbf{(a):} Training loss of the candidate solution during fine-tuning. \textbf{(b):} Relative error of the candidate solution.}%
    \label{fig:3.5}%
\end{figure}

Figure~\ref{fig:3.5} shows the training loss and relative error during the fine-tuning of the candidate function. The solution identified by FEX is $\tilde{u}(t,x) = -0.0000392 + 0.999998 x$, which closely approximates the true solution $u(t,x) = x$. The relative error of $\tilde{u}$ is $4.87 \times 10^{-6}$, which is a significant improvement over the TD-NN results of $\mathcal{O}(10^{-3})$ reported in \cite{lu}.

{
\subsection{The Variance Gamma Model}
To further demonstrate our method's strength, we present another example.  While in low dimension, this problem shows FEX's strength at solving problems beyond toy examples and function approximation.  In addition, we use this problem to verify the proposed integral simplification method given in Section~\ref{sec:evalint}.  Motivated by \cite{Fu03082022}, we seek to solve the Variance Gamma (VG) options pricing problem.  The VG model is a widely used pure-jump L\'{e}vy process that captures the observed features of asset returns that are not well modeled by the classical Black-Scholes framework.  In the context of European option pricing, the VG model leads to a PIDE.  Letting $V(S,t)$ be the value of the option given asset price $S$ and time until option maturity $t$ the PIDE governing the model's behavior can be written as
\begin{align}
    \int_{-\infty}^{\infty} \Big (V(Se^y,t) - V(S,t) - \frac{\partial V}{\partial S}(S,t)S(e^y - 1)\Big ) k(y)dy + \frac{\partial V}{\partial t}(S,t) + (r-q)S\frac{\partial V}{\partial S}(S,t) - rV(S,t) = 0.
    \label{eqn:VGPIDE}
\end{align}

The measure $k(y)$ is the L\'{e}vy density of the VG process, which is given by
\begin{align}
k(y) = \frac{e^{-\lambda_p y}}{\nu y} 1_{y>0} + \frac{e^{-\lambda_n |y|}}{\nu |y|} 1_{y<0},
\label{VGmeasure}
\end{align}
where 
\begin{align}
    \lambda_p = \Big(\frac{\theta^2}{\sigma^4} + \frac{2}{\sigma^2 \nu}\Big)^{\frac{1}{2}} - \frac{\theta}{\sigma^2},
    \label{lambdap}
\end{align}
and
\begin{align}
\lambda_n = \Big(\frac{\theta^2}{\sigma^4} + \frac{2}{\sigma^2 \nu}\Big)^{\frac{1}{2}} + \frac{\theta}{\sigma^2}.
\label{lambdan}
\end{align}
We choose $\theta = -0.4$, $\sigma = 0.4$, $\nu = 0.4$, $r = 0.05$ and $q = 0.02$, all as in \cite{Fu03082022}.  To facilitate training we rescale the problem so that the asset value $S$ lies within $[0,1]$.  Our strike price is likewise scaled, and is then $K = 1/3$.  Time to maturity, $t$ is within $[0,2]$.  Note that $r$ is the risk-free interest rate and $q$ is the dividend rate of the stock. We will model European puts, so our initial condition is 
$V(0,S) = (S-K)^{+}$ where 
$(S-K)^{+} = \max(S-K, 0)$.  Because of this initial condition, another change made is to include ReLU as a unary function that FEX can use to represent the solution - this choice is straightforward given the form of $\max(S-K,0)$, and will allow a better fit while keeping expression size short.  The boundary conditions are given as $V(t, 0) = Ke^{-rt}$ and $V(t,1) = 0$.  Note that we are not making the common change of variables $x = log(S)$.  Rather than use log-price as our variable we simply use raw asset price.  The functional $\mathcal{L}$ \eqref{eqn:leastsquare} now includes an extra term as we have both boundary conditions and an initial condition.  Letting $\mathcal{D}$ once again represent the PIDE, in this case Eqn. \eqref{eqn:VGPIDE}, we have:
\begin{align}
    \mathcal{L}(u) \approx \frac{1}{N}\sum_{i=1}^N|\mathcal{D}(\tilde{u}(t_i,S_i))|^2+ \frac{1}{M}\sum_{j=1}^M | \tilde{u}(t_j, S_{min}) - Ke^{-rt}|^2 + |\tilde{u}(t_j, S_{max})|^2+|\tilde{u}(0, S_j) - (S_j - K)^{+}|^2.
    \label{eqn:VGleastsquare}
\end{align} 
As noted earlier, to evaluate the integral term we employ the same Taylor approximation trick as proposed in Section \ref{sec:evalint}, with the slight difference being that we need to find $\mathbb{E}[{u}(t, xG(x,y))] = \mathbb{E}[{u}(t,xe^y)]$  rather than $\mathbb{E}[{u}(t,x + z)]$. 

To benchmark our results, we will compare FEX's results to those from the Carr-Madan Fast Fourier transform based method \cite{CarrMadan1999}.  It is important to note that an explicit solution to this problem does not exist.  Hence FEX is learning an approximation, which greatly increases the challenge.  

\begin{figure}[!ht]
    \centering
    \includegraphics[scale=.29]{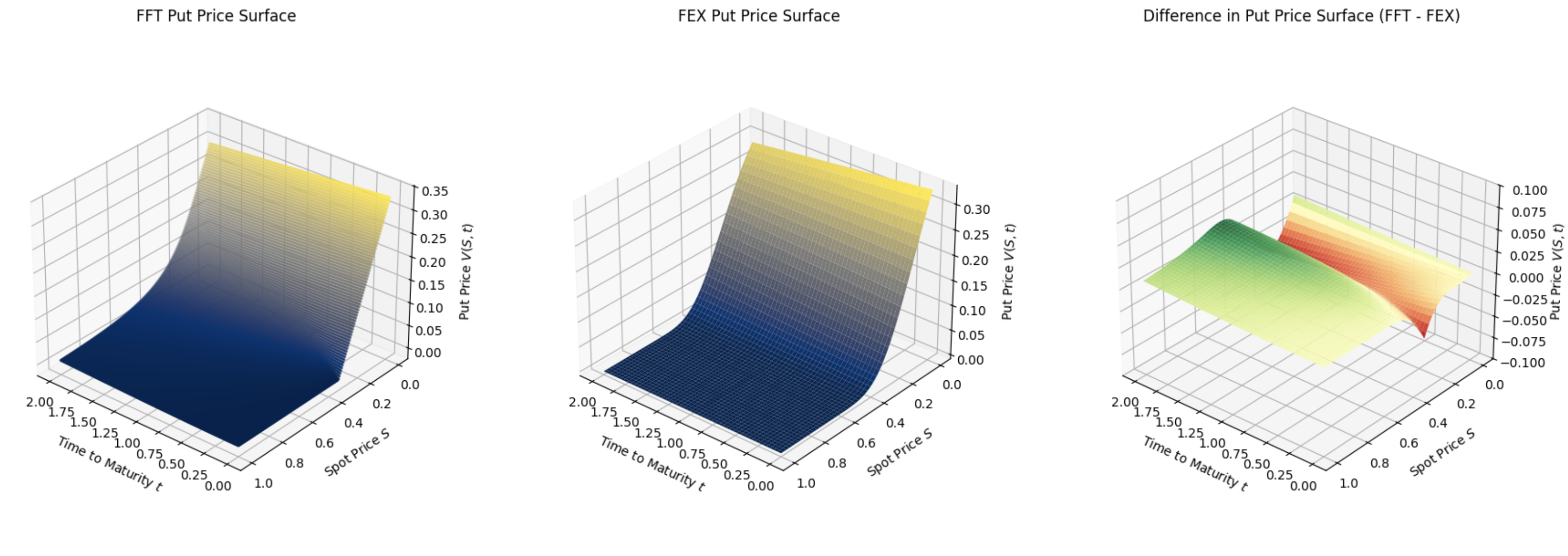}
    \caption{{Plot of the put value surfaces given by the Carr-Madan FFT method \cite{CarrMadan1999}, FEX, and the difference between the two.}}
    \label{fig:VG}
\end{figure}

The true solution, the numerical solution, and their difference is visulaized in Figure \ref{fig:VG}. The average mean square error (measured across ten consecutive runs of FEX) between FEX's solution and the Carr-Madan solution was $2.25\times10^{-4}$, indicating a good fit.  300 iterations of the searching loop were performed using a depth-3 tree (one level deeper than the tree seen on the far right of Figure~\ref{fig:tree}). Due to the lack of an explicit true solution we did not use a stopping criteria based on accuracy; rather, we simply used a fixed number of fine-tune steps (30,000). Importantly, no significant change in accuracy was found when numerical integration (using the trapezoid rule) was performed instead of the simplification discussed in Section \ref{sec:evalint}, but the speedup was very significant.
}
\subsection{High Dimensional PIDE}
\label{sec:highdpide}
Having established the strength of FEX in solving PIDEs and functional approximation in lower dimensions, all newly proposed techniques in FEX, including simplifying the integral and using parameter grouping as introduced in Section~\ref{sec:alg},are adopted to extend the excellent performance of FEX from lower-dimensional cases to this high-dimensional problem.

\textbf{Example 1. } Consider the following PIDE:
\begin{align}
    \begin{cases} 
    \frac{\partial u}{\partial t} + \frac{\epsilon}{2} x \cdot \nabla u(t,x) + \frac{1}{2} Tr(\sigma^2 H(u)), \\
    + \int_{\mathrm{R}^d} (u(t, x+z) - u(t,x) - z \cdot \nabla u(t,x))\nu(dz) = \lambda(\mu^2 + \sigma^2) + \theta^2 + \frac{\epsilon}{2}||x||^2,\\
    u(T, x) = \frac{1}{d}||x||^2.\\
    \end{cases}
    \label{eqn:3.9}
\end{align}
Here, as in \cite{lu}, $\epsilon = 0$, $\theta = 0.3$, Poisson intensity $\lambda = 0.3$, $\mu = 1$, $\sigma^2 = 0.0001$. The true solution is $\frac{1}{d}||x||^2$, for $x \in \mathbb{R}^d$. In this example, various dimensions $d$ are tested and the parameter grouping method introduced in Section \ref{sec:paramgrouping} is used to accelerate coefficient optimization. 


\begin{figure}[!ht]
    \centering
    \subfloat[\centering Relative Error]{{\includegraphics[width=7cm]{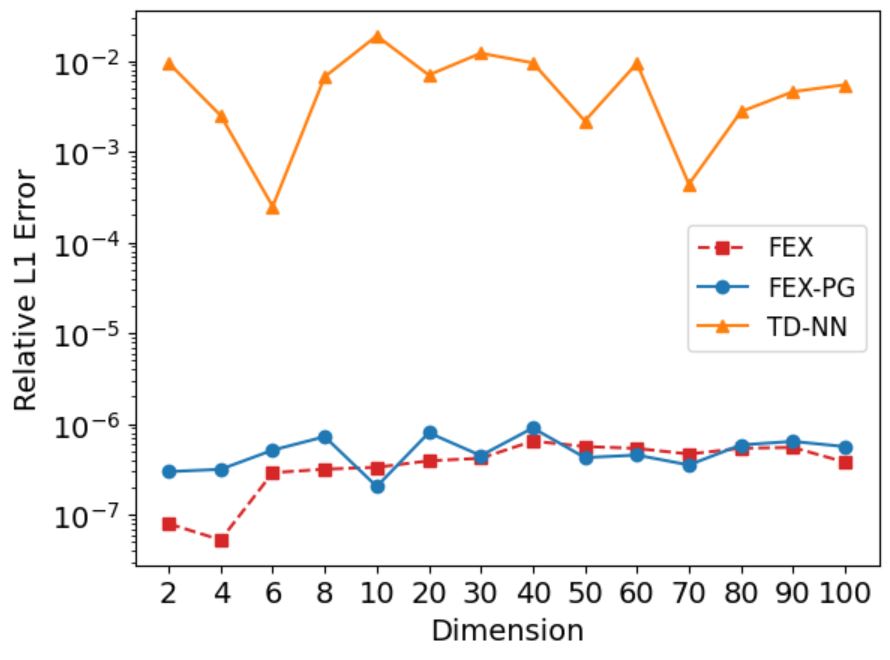}}\label{fig:3.9relerror}}%
    \qquad
    \subfloat[\centering Run Time]{{\includegraphics[width=7.3cm]{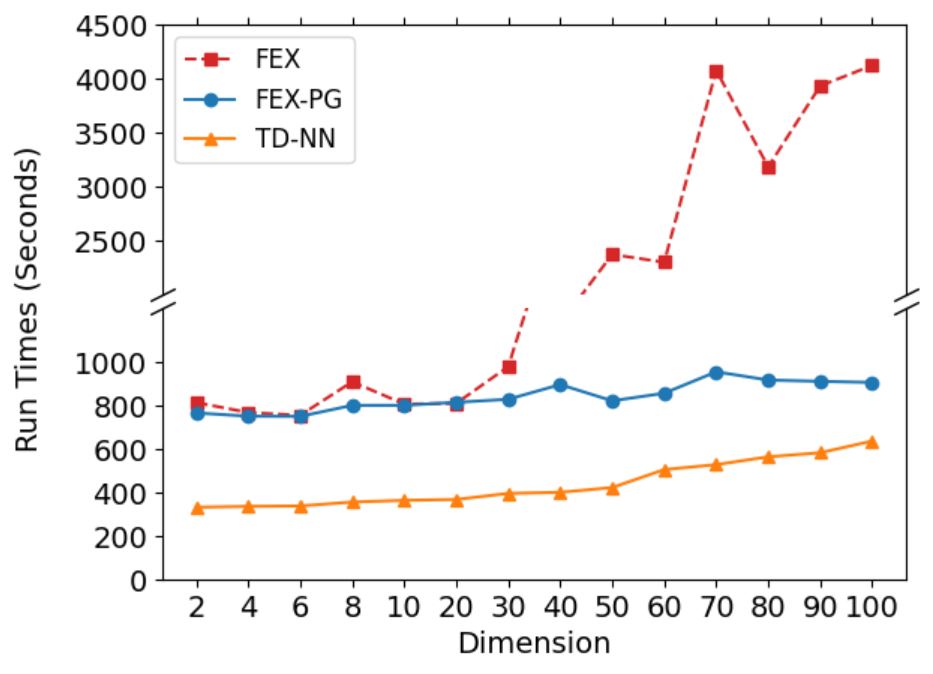}}\label{fig:3.9runtime}}%
    \caption{Accuracy and run time comparison for PIDE~\eqref{eqn:3.9} across different problem dimensions.}%
    \label{fig:3.9}%
\end{figure}

Figure \ref{fig:3.9relerror} (Table~\ref{tab:err1}) shows a comparison of the relative error in the solutions produced by FEX-based methods and TD-NN \cite{lu} across various problem dimensions. Since single-precision floating-point numbers are used, the relative error of the FEX-based methods reaches the expected single-precision accuracy, whereas TD-NN has an error of $\mathcal{O}(10^{-2})$. Notably, FEX-PG maintains low relative error even beyond 100 dimensions, achieving a relative error of $1.46\times 10^{-6}$ at 500 dimensions, far surpassing the accuracy of TD-NN. 


Figure \ref{fig:3.9runtime} (Table~\ref{tab:time1}) compares the average runtime of FEX-based methods with that of TD-NN \cite{lu} across various dimensions. 
{ For all experiments, including these timed high dimensional runs, FEX-PG was run on a single NVIDIA A100 80G.  As seen in Algorithm \ref{alg1}, threshold parameter $\eta = \frac{1}{d}$ where $d$ is the number of space dimensions( i.e. $x \in \mathrm{R}^d$).  50 iterations of the searching loop were performed, and the fine-tuning was stopped when the loss functional $\mathcal{L}$ was below $1.5\times10^{-14}$.  The cut off (of near to double precision epsilon) corresponded well with the found solution having a relative error near to single precision epsilon.  This allowed us to investigate how much time was saved with the parameter grouping step while holding the accuracy of the found solution largely constant, as seen in Figure \ref{fig:3.9relerror} and Figure \ref{fig:3.9runtime} respectively.
} 
While TD-NN demonstrates faster performance, FEX delivers significantly higher accuracy, making the trade-off in runtime justifiable. By introducing parameter grouping (FEX-PG, ~\ref{alg1}), the computational cost is reduced, especially in higher dimensions. This improvement is driven by two key factors: first, it becomes easier to identify a good candidate solution during the search loop, particularly in very high dimensions (80+); second, parameter grouping allows us to optimize far fewer parameters, significantly reducing the number of optimization iterations required to achieve similar accuracy. In this case, single precision accuracy within approximately 2,000 iterations of Adam is achieved, compared to as many as 20,000 iterations without parameter grouping.  

\begin{figure}[!ht]
    \centering
    \includegraphics[scale=.35]{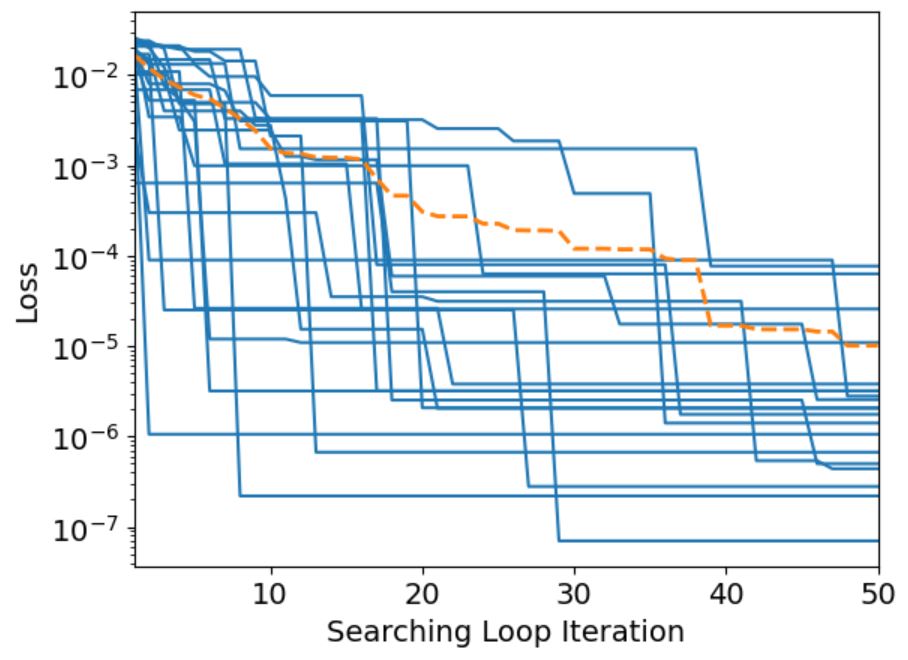}
    \caption{Loss during searching phase of FEX-PG, while solving \eqref{eqn:3.9} in 100 dimensions. Each blue line represents an independent experiment.  The lines plot the loss of the best candidate in the pool, as consecutive iterations of the searching loop are performed.  20 such trajectories are pictured, with the average shown in dashed orange.}
    \label{fig:searchingloss}
\end{figure}

Figure \ref{fig:searchingloss} illustrates the process of the searching loop used to solve the CO problem of finding a correct sequence of operators.  As successive iterations are performed, the pool of functions includes better and better candidates.  In practice, it is observed that a loss on the order of $10^{-4}$ or $10^{-5}$ is easily good enough to indicate a good candidate function.


\begin{figure}[!ht]
\centering
    \subfloat[]{\includegraphics[width=0.45\linewidth]{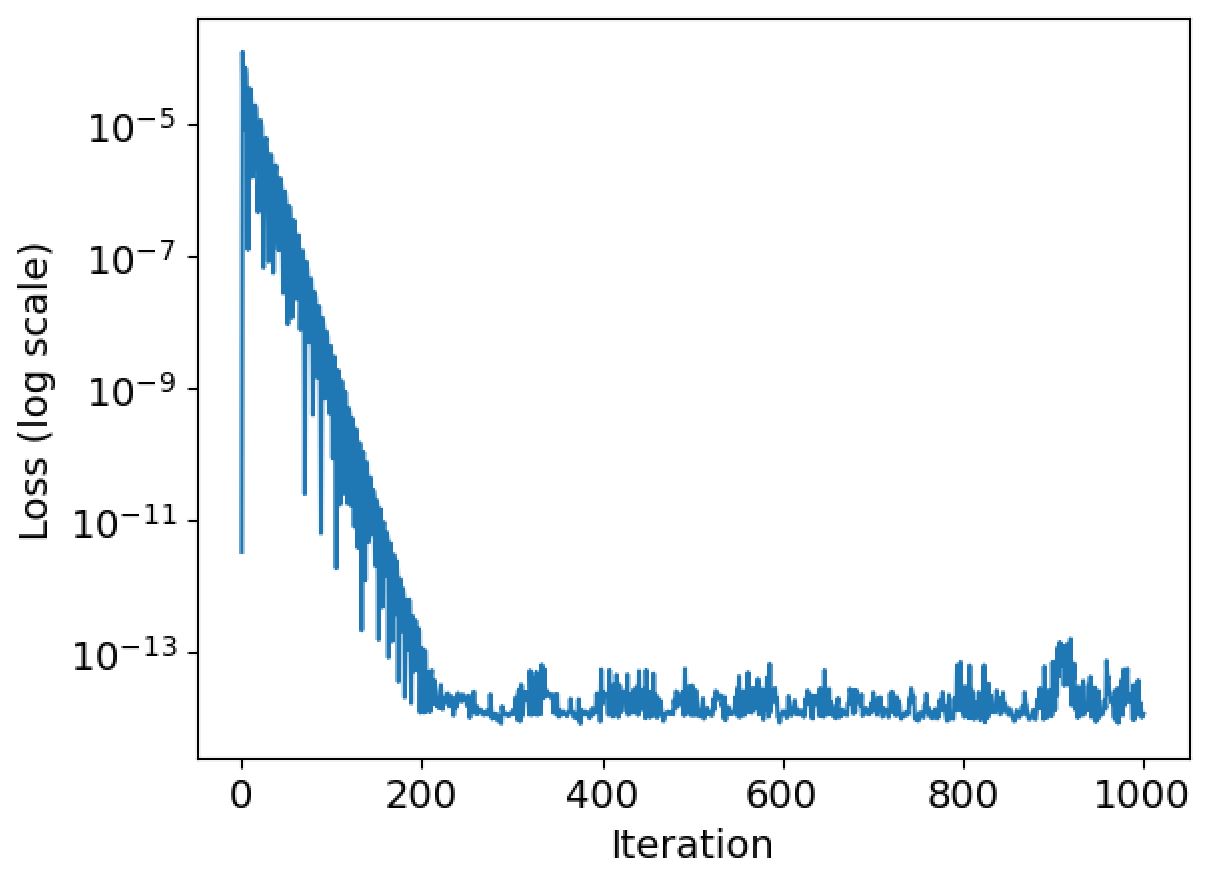}}
\hfil
    \subfloat[]{\includegraphics[width=0.45\linewidth]{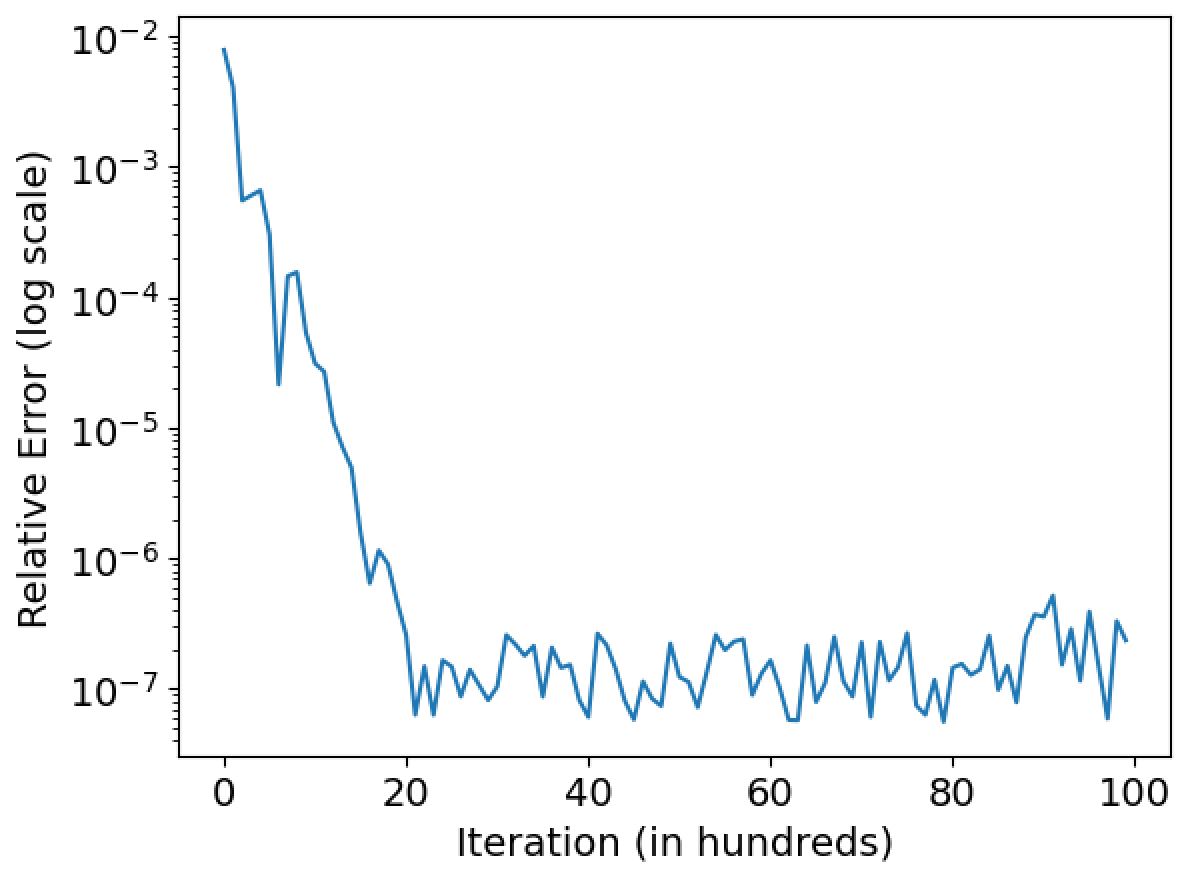}}

    \subfloat[]{\includegraphics[width=0.45\linewidth]{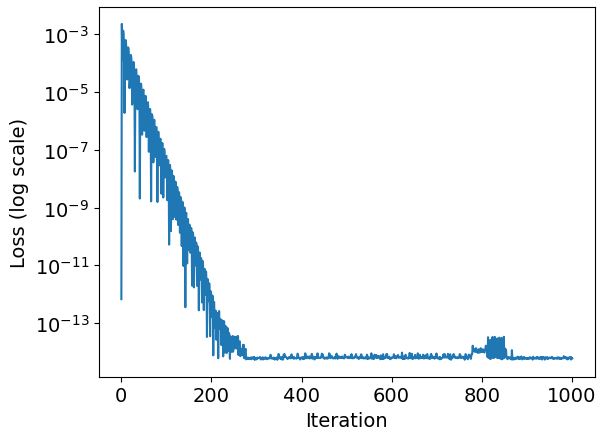}}
\hfil
    \subfloat[]{\includegraphics[width=0.45\linewidth]{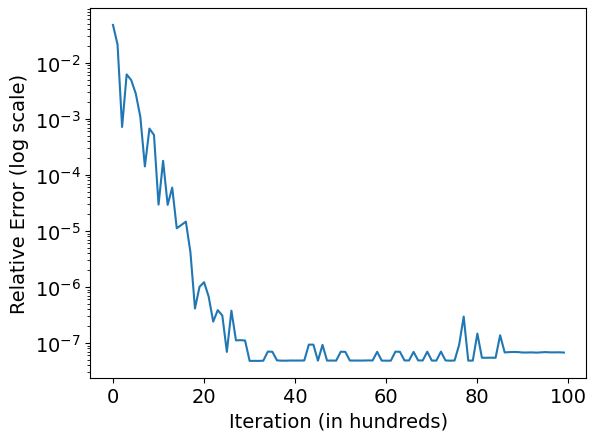}}

    \subfloat[]{\includegraphics[width=0.45\linewidth]{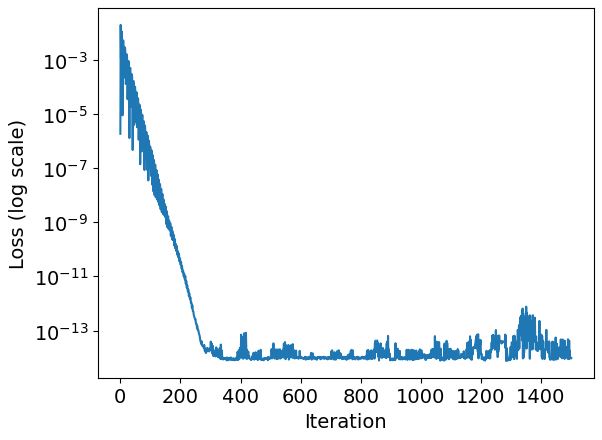}}
\hfil
    \subfloat[]{\includegraphics[width=0.45\linewidth]{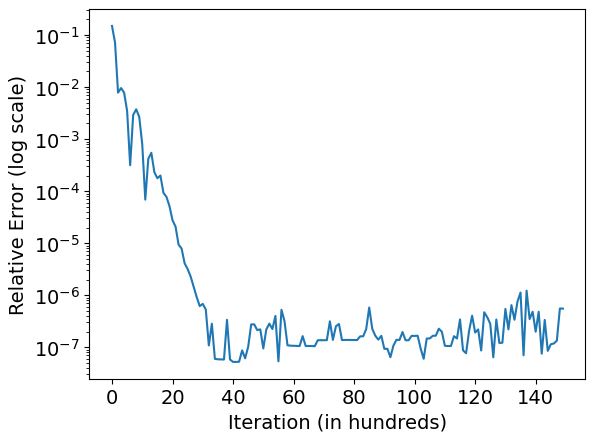}}
\caption{Optimization profile for PIDE~\eqref{eqn:3.9} of 10 (Row 1), 50 (Row 2) and 100 (Row 3) dimensions. \textbf{Col 1:} Training loss of the candidate solution during fine-tuning. \textbf{Col 2:} Relative error of the candidate solution during fine-tuning.}
    \label{fig:my figure}
    \end{figure}



\begin{table}[ht]
\begin{center}
\begin{tabular}{||c | c c c c c c c||} 
 \hline
 Dimension & 2 & 4 & 6 & 8 & 10 & 20 & 30 \\ [0.5ex] 
 \hline\hline
 FEX-PG & 2.99e-7 & 3.17e-7 & 5.16e-7 & 7.26e-7 & 2.05e-7 & 8.02e-7 & 4.49e-7 \\ 
 \hline
 TD-NN \cite{lu} & 0.00954 & 0.00251 & 0.00025 & 0.00671 & 0.01895  & 0.00702 & 0.01221 \\
 \hline 
\end{tabular}
\end{center}
\begin{center}
\begin{tabular}{||c | c c c c c c c||}  
 \hline
 Dimension & 40 & 50 & 60 & 70 & 80 & 90 & 100\\ [0.5ex] 
 \hline\hline
 FEX-PG & 9.05e-7 & 4.27e-7 & 4.55e-7 & 3.54e-7 & 5.89e-7 & 6.44e-7 &  5.64e-7 \\ 
 \hline
 TD-NN \cite{lu} &  0.00956 & 0.00219 &  0.00944 & 0.00044 & 0.00277 & 0.00460 & 0.00548 \\
 \hline 
\end{tabular}
\end{center}
\caption{The comparision of the FEX-PG and TD-NN \cite{lu} in terms of relative errors.}
\label{tab:err1}
\end{table}

\begin{table}[ht]
\begin{center}
\begin{tabular}{||c | c c c c c c c||} 
 \hline
 Dimension & 2 & 4 & 6 & 8 & 10 & 20 & 30 \\ [0.5ex] 
 \hline\hline
 FEX-PG & 768 & 753 & 752 & 803 & 803 & 817 & 831 \\ 
 \hline
 TD-NN \cite{lu} & 332 & 336 & 338 & 356 & 364  & 368 & 396 \\
 \hline 
\end{tabular}
\end{center}

\begin{center}
\begin{tabular}{||c | c c c c c c c||} 
 \hline
 Dimension & 40 & 50 & 60 & 70 & 80 & 90 & 100 \\ [0.5ex] 
 \hline\hline
 FEX-PG & 899 & 824 & 859 & 958 & 920 & 914 & 909 \\ 
 \hline
 TD-NN \cite{lu} & 401 & 423 & 506 & 529 & 565  & 584 & 638 \\
 \hline 
\end{tabular}
\end{center}
\caption{The comparision of the FEX-PG and TD-NN \cite{lu} in terms of the computation times.}
\label{tab:time1}
\end{table}

Due to the use of the integral estimation~\eqref{eqn:3.9}, the error introduced by the integral estimation increases with the variance $\sigma$. To assess the performance and asses the robustness of FEX-PG under higher variance scenarios, larger values of $\sigma$ are tested. Table~\ref{tab:var} compares the relative error of the solution obtained by FEX-PG and TD-NN for PIDE~\eqref{eqn:3.9} in $d=100$ dimensions across various levels of variance. Even with a variance of $1$, FEX-PG performs exceptionally well, maintaining a relative error on the order of $10^{-7}$, whereas the error of TD-NN increases as the variance grows. 
{
This good behavior of FEX is also demonstrated by the earlier example of the VG model, where the variance of $k(y)$, the L\'{e}vy density of the VG process is $0.064$.
}
It is important to note that, since the integral term is approximated using a Taylor series with respect to the random variable $z$, the use 
{
of heavy-tailed distributions
}
for $z$ may cause instability, due to the the need of finite moments (without this, the remainder of the Taylor series may not converge). While this limitation may inspire future research, it remains unexplored for the time being.

\begin{table}[H]
\centering
\begin{tabular}{||c | c c c c ||} 
 \hline
 Variance & 0.001 & 0.01 & 0.1 & 1.0 \\ [0.5ex] 
 \hline\hline
 FEX-PG & 4.07e-7 & 4.71e-7 & 5.83e-7 & 3.05e-7\\ 
 \hline
 TD-NN \cite{emailLiweiLu} & 0.00446 & 0.00639 & 0.02716 & 0.03821 \\
 \hline 
\end{tabular}
\caption{The relative errors of the solution obtained by FEX-PG and TD-NN for PIDE~\eqref{eqn:3.9} in $d=100$ dimensions across various levels of variance. }
\label{tab:var}
\end{table}

\textbf{Example 2.} Finally we arrive at the last example:
\begin{align}
    \begin{cases} 
    \frac{\partial u}{\partial t} + \frac{\epsilon}{2} ||x|| x \cdot \nabla u(t,x) + \frac{1}{2} Tr(\boldsymbol{\sigma} \boldsymbol{\sigma}^T H(u)) \\
    + \int_{\mathrm{R}^d} (u(t, x+z) - u(t,x) - z \cdot \nabla u(t,x))\nu(dz) = \lambda(\mu^2 + \sigma^2) + \frac{2d-2}{d} \theta^2 + \frac{\epsilon}{d}||x||^3,\\
    u(T, x) = \frac{1}{d}||x||^2.\\
    \end{cases}
    \label{eqn:3.10}
\end{align}
Once again we use identical constants to \cite{lu}, where $\epsilon = 0.05$, $\theta = 0.2$, $\lambda = 0.3$, $\sigma = 0.0001$ and the matrix $\boldsymbol{\sigma}$ is given as:
\begin{align*}
\boldsymbol{\sigma} = \theta 
\begin{bmatrix} 
1 & 0 & 0 & 0 & \cdots & 0 \\
1 & 1 & 0 & 0 & \cdots & 0 \\
0 & 1 & 1 & 0 & \cdots & 0 \\
0 & 0 & 1 & 1 & \cdots & 0 \\
\vdots & \vdots & \vdots & \vdots & \ddots & 0 \\
0 & 0 & 0 & \cdots & 1 & 1 \\
\end{bmatrix}.    
\end{align*}

\begin{figure}[ht]
\centering
    \subfloat[]{\includegraphics[width=0.45\linewidth]{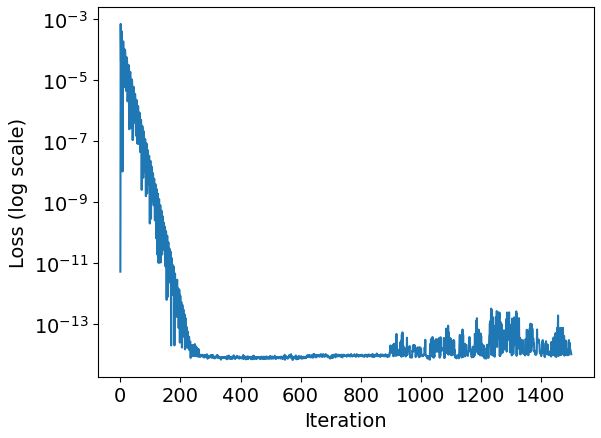}}
\hfil
    \subfloat[]{\includegraphics[width=0.45\linewidth]{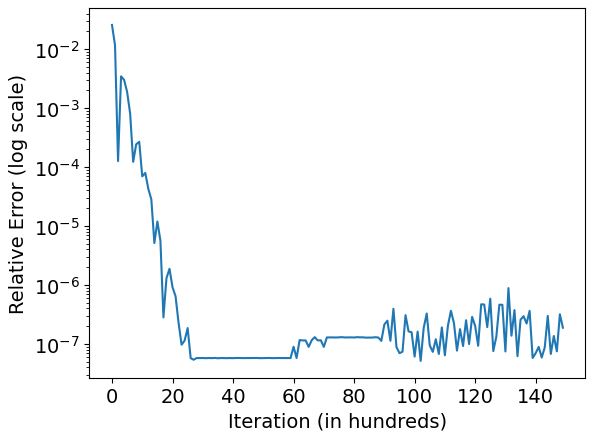}}

    \subfloat[]{\includegraphics[width=0.45\linewidth]{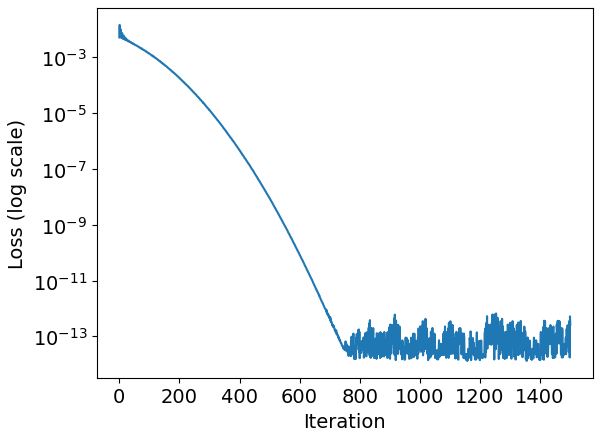}}
\hfil
    \subfloat[]{\includegraphics[width=0.45\linewidth]{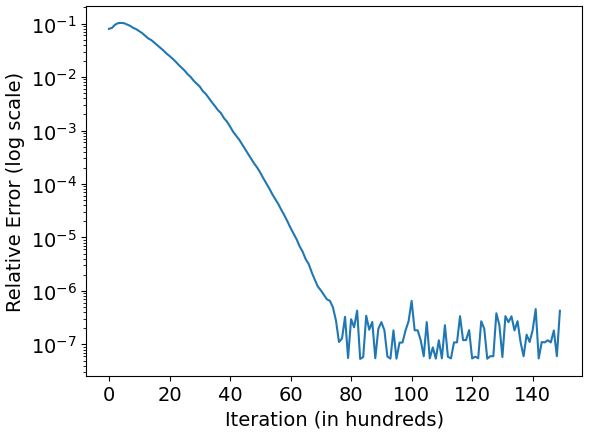}}
\caption{Optimization profile for PIDE~\eqref{eqn:3.10} of 25 (Row 1) and 100 (Row 2) dimensions. \textbf{Col 1:} Training loss of the candidate solution during fine-tuning. \textbf{Col 2:} Relative error of the candidate solution.}
    \label{fig:3.10}
    \end{figure}

Table~\ref{tab:3.10} compares our results with those presented in \cite{lu}. Once again, accuracy of several orders of magnitude higher is achieved. It is worth noting that the FEX-PG method converged to these results in roughly the same time as for Eqn.~\eqref{eqn:3.9}, indicating that the increased complexity of the problem did not notably impact the learning process—neither during the search phase for the controller nor during the fine-tuning of the weights. Figure~\ref{fig:3.10} demonstrates rapid convergence during fine-tuning, achieving single precision accuracy within just 2,000 iterations of Adam.
\begin{table}[H]
    \centering
    \begin{tabular}{||c | c c c c ||} 
 \hline
 Dimension & 25 & 50 & 75 & 100 \\ [0.5ex] 
 \hline\hline
 FEX-PG & 6.97e-8 & 1.67e-7 & 9.51e-8 & 1.60e-7 \\ 
 \hline
 TD-NN \cite{lu} & 0.00743 & 0.01910 & 0.02412 & 0.02387 \\
 \hline 
\end{tabular}
    \caption{The relative errors of the solution obtained by FEX and TD-NN for PIDE~\eqref{eqn:3.9} across various dimension.}
    \label{tab:3.10}
\end{table}

\subsection{Conclusion and Discussion}

This paper presents a numerical approach to solving the PIDE using a new finite expression method, the FEX-PG. While the standard finite expression method works well for solving low-dimensional PIDEs, its efficiency decreases significantly during the search loop and parameter optimization as the dimensionality increases. To address this, two improvements are proposed: first, simplifying the estimation of the integral term, and second, grouping trainable parameters in the candidate solution to reduce the total number of parameters that need to be optimized. Through extensive experiments on problems with various dimensions (up to 500), the method's effectiveness in accurately identifying the true solution and achieving machine-level precision is demonstrated. In contrast, the baseline neural network approach struggles with higher dimensions and variance, a key parameter in the PIDE.

As a preliminary demonstration, the focus of the high dimensional experiments has been on the specific case where $G(x,z)=z$ in the integral term~\eqref{eqn:integral}. This serves as a proof of concept to validate the effectiveness of our proposed method. In future work, this approach may be extended by exploring more general forms of $G(x,z)$, which will allow for greater flexibility and applicability of the FEX-PG to a wider range of problems. The corresponding estimation techniques for the integral term~\eqref{eqn:integral} for these more complex choices of $G(x,z)$ will be developed.


\section*{Acknowledgement}
S. W. was partially supported by the U.S. Department of Energy, Office of Science, Office of Advanced Scientific Computing Research SciDAC program under Contract No. DE-AC02-05CH11231. Y. H. was partially supported by the US National Science Foundation under awards DMS-2244988, DMS-2206333, the Office of Naval Research Award N00014-23-1-2007, and the DARPA D24AP00325-00. Approved for public release; distribution is unlimited.
\bibliographystyle{plain} 
\bibliography{ref.bib}

\begin{thebibliography}{10}

\bibitem{Applebaum}
D.~Applebaum.
\newblock {\em L\'{e}vy Processes and Stochastic Calculus}.
\newblock Cambridge University Press, 2009.

\bibitem{bellman}
Richard Bellman and Robert Kalaba.
\newblock A mathematical theory of adaptive control processes.
\newblock {\em Proceedings of the National Academy of Sciences}, 45(8):1288--1290, 1959.

\bibitem{CarrMadan1999}
Peter Carr and Dilip~B. Madan.
\newblock Option valuation using the fast fourier transform.
\newblock {\em Journal of Computational Finance}, 2(4):61--73, 1999.

\bibitem{chizat}
Lenaic Chizat and Francis Bach.
\newblock Implicit bias of gradient descent for wide two-layer neural networks trained with the logistic loss, 2020.

\bibitem{cruz}
Jose Cruz, Maria Grossinho, Daniel Sevcovic, and Cyril~Izuchukwu Udeani.
\newblock Linear and nonlinear partial integro-differential equations arising from finance, 2022.

\bibitem{cuomo}
Salvatore Cuomo, Vincenzo~Schiano Di~Cola, Fabio Giampaolo, Gianluigi Rozza, Maziar Raissi, and Francesco Piccialli.
\newblock Scientific machine learning through physics--informed neural networks: Where we are and what's next.
\newblock {\em Journal of Scientific Computing}, 92(3):88, Jul 2022.

\bibitem{E2017}
Weinan E, Jiequn Han, and Arnulf Jentzen.
\newblock Deep learning-based numerical methods for high-dimensional parabolic partial differential equations and backward stochastic differential equations.
\newblock {\em Communications in Mathematics and Statistics}, 5(4):349--380, Dec 2017.

\bibitem{e2017deepritzmethoddeep}
Weinan E and Bing Yu.
\newblock The deep ritz method: A deep learning-based numerical algorithm for solving variational problems, 2017.

\bibitem{fletcher2013practical}
Roger Fletcher.
\newblock {\em Practical methods of optimization}.
\newblock John Wiley \& Sons, 2013.

\bibitem{Fu03082022}
Weilong Fu and Ali~Hirsa and.
\newblock An unsupervised deep learning approach to solving partial integro-differential equations.
\newblock {\em Quantitative Finance}, 22(8):1481--1494, 2022.

\bibitem{Grohs_2023}
Philipp Grohs, Fabian Hornung, Arnulf Jentzen, and Philippe von Wurstemberger.
\newblock A proof that artificial neural networks overcome the curse of dimensionality in the numerical approximation of black–scholes partial differential equations.
\newblock {\em Memoirs of the American Mathematical Society}, 284(1410), April 2023.

\bibitem{doi:10.1073/pnas.1718942115}
Jiequn Han, Arnulf Jentzen, and Weinan E.
\newblock Solving high-dimensional partial differential equations using deep learning.
\newblock {\em Proceedings of the National Academy of Sciences}, 115(34):8505--8510, 2018.

\bibitem{HAN2020109792}
Jiequn Han, Jianfeng Lu, and Mo~Zhou.
\newblock Solving high-dimensional eigenvalue problems using deep neural networks: A diffusion monte carlo like approach.
\newblock {\em Journal of Computational Physics}, 423:109792, 2020.

\bibitem{jiang}
Zhongyi Jiang, Chunmei Wang, and Haizhao Yang.
\newblock Finite expression methods for discovering physical laws from data.

\bibitem{Karniadakis2021}
George~Em Karniadakis, Ioannis~G. Kevrekidis, Lu~Lu, Paris Perdikaris, Sifan Wang, and Liu Yang.
\newblock Physics-informed machine learning.
\newblock {\em Nature Reviews Physics}, 3(6):422--440, Jun 2021.

\bibitem{KHOO_LU_YING_2021}
Yuehaw Khoo, Jianfeng Lu, and Lexing Ying.
\newblock Solving parametric pde problems with artificial neural networks.
\newblock {\em European Journal of Applied Mathematics}, 32(3):421–435, 2021.

\bibitem{kitora}
Shuji Kitora, Souma Jinno, Hiroshi Toki, and Masayuki Abe.
\newblock Relation of integro-partial differential equations with delay effect based on the maxwell equations to the heaviside and pocklington equations.
\newblock {\em IEEE Transactions on Electromagnetic Compatibility}, 63(4):1223--1230, 2021.

\bibitem{lapa}
Eugenio Lapa, Adrian LL, Willy BarahonaM, and Gabriel V.
\newblock Integro-differential equation of p-kirchhoff type with no-flux boundary condition and nonlocal source term.
\newblock {\em Journal of Applied Mathematics and Mechanics}, 2:23--30, 01 2015.

\bibitem{Li_2022}
Haoya Li and Lexing Ying.
\newblock A semigroup method for high dimensional elliptic pdes and eigenvalue problems based on neural networks.
\newblock {\em Journal of Computational Physics}, 453:110939, March 2022.

\bibitem{liang2022stiffness}
Senwei Liang, Zhongzhan Huang, and Hong Zhang.
\newblock Stiffness-aware neural network for learning hamiltonian systems.
\newblock In {\em International Conference on Learning Representations}, 2022.

\bibitem{liang2024solving}
Senwei Liang, Shixiao~W Jiang, John Harlim, and Haizhao Yang.
\newblock Solving pdes on unknown manifolds with machine learning.
\newblock {\em Applied and Computational Harmonic Analysis}, 71:101652, 2024.

\bibitem{liang}
Senwei Liang and Haizhao Yang.
\newblock Finite expression method for solving high-dimensional partial differential equations.
\newblock {\em ArXiv, abs/2206.10121}, 2023.

\bibitem{Ming_2021}
Yulei Liao and Pingbing Ming.
\newblock Deep nitsche method: Deep ritz method with essential boundary conditions.
\newblock {\em Communications in Computational Physics}, 29(5):1365–1384, June 2021.

\bibitem{emailLiweiLu}
Liwei Lu, Hailong Guo, Xu~Yang, and Yi~Zhu.
\newblock private communication, May 2024.
\newblock High variance results taken from private communication with corresponding author.

\bibitem{lu}
Liwei Lu, Hailong Guo, Xu~Yang, and Yi~Zhu.
\newblock Temporal difference learning for high-dimensional pides with jumps.
\newblock {\em SIAM Journal on Scientific Computing}, 46(4):C349--C368, 2024.

\bibitem{lu2021machinelearningellipticpdes}
Yiping Lu, Haoxuan Chen, Jianfeng Lu, Lexing Ying, and Jose Blanchet.
\newblock Machine learning for elliptic pdes: Fast rate generalization bound, neural scaling law and minimax optimality, 2021.

\bibitem{NEURIPS2020_2000f632}
Yulong Lu and Jianfeng Lu.
\newblock A universal approximation theorem of deep neural networks for expressing probability distributions.
\newblock In H.~Larochelle, M.~Ranzato, R.~Hadsell, M.F. Balcan, and H.~Lin, editors, {\em Advances in Neural Information Processing Systems}, volume~33, pages 3094--3105. Curran Associates, Inc., 2020.

\bibitem{pmlr-v134-lu21a}
Yulong Lu, Jianfeng Lu, and Min Wang.
\newblock A priori generalization analysis of the deep ritz method for solving high dimensional elliptic partial differential equations.
\newblock In Mikhail Belkin and Samory Kpotufe, editors, {\em Proceedings of Thirty Fourth Conference on Learning Theory}, volume 134 of {\em Proceedings of Machine Learning Research}, pages 3196--3241. PMLR, 15--19 Aug 2021.

\bibitem{müllner2011modernhierarchicalagglomerativeclustering}
Daniel Müllner.
\newblock Modern hierarchical, agglomerative clustering algorithms, 2011.

\bibitem{sym12081379}
Lukáš Novák and Drahomír Novák.
\newblock On taylor series expansion for statistical moments of functions of correlated random variables.
\newblock {\em Symmetry}, 12(8), 2020.

\bibitem{OLSON1991309}
Jane~M. Olson and Lisa~A. Weissfeld.
\newblock Approximation of certain multivariate integrals.
\newblock {\em Statistics {\&} Probability Letters}, 11(4):309--317, 1991.

\bibitem{petersen2021deep}
Brenden~K Petersen, Mikel~Landajuela Larma, Terrell~N. Mundhenk, Claudio~Prata Santiago, Soo~Kyung Kim, and Joanne~Taery Kim.
\newblock Deep symbolic regression: Recovering mathematical expressions from data via risk-seeking policy gradients.
\newblock In {\em International Conference on Learning Representations}, 2021.

\bibitem{raissi}
M.~Raissi, P.~Perdikaris, and G.E. Karniadakis.
\newblock Physics-informed neural networks: A deep learning framework for solving forward and inverse problems involving nonlinear partial differential equations.
\newblock {\em Journal of Computational Physics}, 378:686--707, 2019.

\bibitem{Shen_2021}
Zuowei Shen, Haizhao Yang, and Shijun Zhang.
\newblock Deep network with approximation error being reciprocal of width to power of square root of depth.
\newblock {\em Neural Computation}, 33(4):1005–1036, 2021.

\bibitem{Shen_2021_3layers}
Zuowei Shen, Haizhao Yang, and Shijun Zhang.
\newblock Neural network approximation: Three hidden layers are enough.
\newblock {\em Neural Networks}, 141:160–173, September 2021.

\bibitem{shen2022deepnetworkapproximationachieving}
Zuowei Shen, Haizhao Yang, and Shijun Zhang.
\newblock Deep network approximation: Achieving arbitrary accuracy with fixed number of neurons, 2022.

\bibitem{Sirignano_2018}
Justin Sirignano and Konstantinos Spiliopoulos.
\newblock Dgm: A deep learning algorithm for solving partial differential equations.
\newblock {\em Journal of Computational Physics}, 375:1339–1364, December 2018.

\bibitem{song}
Zezheng Song, Maria~K. Cameron, and Haizhao Yang.
\newblock A finite expression method for solving high-dimensional committor problems.
\newblock {\em SIAM Journal of Scientific Computing}, 2024.

\bibitem{song2024finiteexpressionmethodlearning}
Zezheng Song, Chunmei Wang, and Haizhao Yang.
\newblock Finite expression method for learning dynamics on complex networks.
\newblock {\em arxiv:2401.03092}, 2024.

\bibitem{tang}
Tao Tang.
\newblock A finite difference scheme for partial integro-differential equations with a weakly singular kernel.
\newblock {\em Applied Numerical Mathematics}, 11(4):309--319, 1993.

\bibitem{vershynin2011introductionnonasymptoticanalysisrandom}
Roman Vershynin.
\newblock Introduction to the non-asymptotic analysis of random matrices, 2011.

\bibitem{2020SciPy-NMeth}
Pauli Virtanen, Ralf Gommers, Travis~E. Oliphant, Matt Haberland, Tyler Reddy, David Cournapeau, Evgeni Burovski, Pearu Peterson, Warren Weckesser, Jonathan Bright, St{\'e}fan~J. {van der Walt}, Matthew Brett, Joshua Wilson, K.~Jarrod Millman, Nikolay Mayorov, Andrew R.~J. Nelson, Eric Jones, Robert Kern, Eric Larson, C~J Carey, {\.I}lhan Polat, Yu~Feng, Eric~W. Moore, Jake {VanderPlas}, Denis Laxalde, Josef Perktold, Robert Cimrman, Ian Henriksen, E.~A. Quintero, Charles~R. Harris, Anne~M. Archibald, Ant{\^o}nio~H. Ribeiro, Fabian Pedregosa, Paul {van Mulbregt}, and {SciPy 1.0 Contributors}.
\newblock {{SciPy} 1.0: Fundamental Algorithms for Scientific Computing in Python}.
\newblock {\em Nature Methods}, 17:261--272, 2020.

\bibitem{wilson}
H.~R. Wilson and J.~D. Cowan.
\newblock A mathematical theory of the functional dynamics of cortical and thalamic nervous tissue.
\newblock {\em Kybernetik}, 13(2):55--80, Sep 1973.

\end{thebibliography}
\end{document}